\documentclass[a4paper,11pt]{article}
\usepackage[dvipdf]{graphicx}
\usepackage[french,english]{babel}
\usepackage{csquotes}
\usepackage{latexsym}
\usepackage{amsmath,amssymb,amsfonts,amsthm}
\usepackage{mathabx}
\usepackage{xcolor}
\usepackage{mathrsfs}
\usepackage{bbm}
\usepackage{float}
\usepackage{tikz}
\usepackage{dsfont}

\usepackage[top=3cm, bottom=3cm, left=2cm, right=2cm]{geometry}

\usepackage{enumerate}

\definecolor{db}{RGB}{0, 0, 130}
\usepackage[colorlinks=true,citecolor=red,linkcolor=db,urlcolor=blue,pdfstartview=FitH]{hyperref}

\usepackage{pdfsync}

\definecolor{rp}{rgb}{0.25, 0, 0.75}
\definecolor{dg}{rgb}{0, 0.5, 0}

\newcommand{\1}{\mathds{1}}
\newcommand{\R}{\mathbb{R}}

\newcommand{\Z}{\mathbb{Z}}
\newcommand{\N}{\mathbb{N}}
\newcommand{\EE}{\mathbb{E}}

\newcommand{\PP}{\mathbb{P}}

\newcommand{\bqn}{\begin{equation}}
\newcommand{\eqn}{\end{equation}}
\newcommand{\bqne}{\begin{equation*}}
\newcommand{\eqne}{\end{equation*}}

\newcommand{\dsum}{\displaystyle \sum}
\newcommand{\apar}[1]{\left(#1\right)}

\newcommand{\ens}[2]{\{#1,\, \cdots, \, #2\}}
\newcommand{\tuple}[2]{\left(#1,\, \cdots, \, #2\right)}
\newcommand{\df}{\mathrm{d}}

\makeatletter
\newcommand{\customlabel}[2]{%
   \protected@write \@auxout {}{\string\newlabel {#1}{{#2}{\thepage}{#2}{#1}{}}}%
   \hypertarget{#1}{#2\hspace{-0.14cm}}
}
\makeatother

\newtheorem{definition}{Definition}[section]
\newtheorem{theorem}[definition]{Theorem}

\newtheorem{corollary}[definition]{Corollary}

\newtheorem{lemma}[definition]{Lemma}

\newtheorem{proposition}[definition]{Proposition}
\newtheorem{remark}[definition]{Remark}
\newtheorem{fact}[definition]{Fact}

\title{Lengths of closed geodesics on random non-compact hyperbolic surfaces}
\author{Malo Gatignol and Joe Thomas}
\date{\today}

\begin{document}
\maketitle
\begin{abstract}
We study the asymptotic distribution of the number of primitive, closed geodesics on Weil-Petersson random hyperbolic surfaces with genus $g$ and $n(g)$ cusps, under various growth rates of $n(g)$ as $g\to\infty$. We prove that the limiting distribution is Poisson with parameter depending explicitly on the limit of $n(g)^2/g$. This work is inspired by that of Mirzakhani and
Petri's \cite{MirzakhaniPetri2021} on closed surfaces and Hide and Thomas \cite{HideThomas2024} when $g$ is fixed. Our proof relies on first order expansions for Weil-Petersson moduli space volumes in the $n(g)/g\to 0$ regime. We also obtain some new existence results for small Laplacian eigenvalues.
\end{abstract}
\tableofcontents
\section{Introduction}

In this paper, we study the distribution of the lengths of primitive closed geodesics on random non-compact hyperbolic surfaces of large area for a variety of different genus versus number of cusp regimes. Our results are inspired by the work of Mirzakhani and Petri \cite{MirzakhaniPetri2021} for closed surfaces and Hide and the second named author \cite{HideThomas2024} for genus zero surfaces with many cusps (see Section \ref{sec: literature} for a broader review of the related literature).

We focus on the Weil-Petersson random model of hyperbolic surfaces which arises from the Weil-Petersson probability measure $\mathbb{P}_{g,n}$ on the moduli space $\mathcal{M}_{g,n}$ consisting of complete, oriented hyperbolic surfaces of genus $g$ with $n$ cusps identified up to isometry (see Section \ref{sec: background} for more details). Let $0\leq a <b$ be real numbers; we are interested in the random variables
    $$N_{g,n,[a,b]}:(\mathcal{M}_{g,n},\mathbb{P}_{g,n})\to \Z_{\geq 0},$$
that count the number of primitive and unoriented closed geodesics on a surface $X\in\mathcal{M}_{g,n}$. 

Akin to counting cycles on random regular graphs, and as observed in the works \cite{MirzakhaniPetri2021,HideThomas2024} for $n=0, g\to\infty$ and $g=0, n\to\infty$, the asymptotic distribution of $N_{g,n,[a,b]}$ in suitable windows $[a,b]$ and regimes of $g+n\to\infty$ is expected to be Poisson. However, the mean of the limiting distribution varies depending upon the growth of $n$ compared to $g$. In this article we will focus on three such regimes and consider $n=n(g)$ to vary as a function of $g\to\infty$.

We will actually prove more than just the limiting distribution of the $N_{g,n,[a,b]}$. Indeed, we will look at the limiting distribution of vectors of these random variables which then allows us to deduce Poisson point process behaviour in the surface length spectra. 

\subsection{Main results}

The first setting that we consider is when $n(g) \ll \sqrt{g}$. Here, the growth rate of the number of cusps is only moderate and not strong enough to affect the limiting behaviour of the random variables from the closed $n=0$ case. The leading order contribution to the moments of the $N_{g,n,[a,b]}$ in this case is given by simple closed geodesics that do not separate the surface when cut along. We hence recover and extend the same asymptotic behaviour as Mirzakhani and Petri in \cite[Theorem 4.1]{MirzakhaniPetri2021}.

\begin{theorem}
   \label{th: stat de poisson dans le cas n^2/g -> 0}
Let $k \in \N$ and
\[a_1 < b_1 < a_2 < b_2 < \ldots < a_k < b_k.\] Assume that
\[\lim_{g\to \infty}\frac{n(g)^2}{g} = 0,\] then
\begin{equation*}
   \tuple{N_{g,n(g), \left[a_1,b_1\right]}}{N_{g,n(g), \left[a_k,b_k\right]}} \xrightarrow[g\to +\infty]{(d)} \tuple{N_{[a_1, b_1]}}{N_{[a_k,b_k]}},
\end{equation*}
where $N_{[a_i, b_i]}$ are pairwise independent random variables following a Poisson distribution with parameter 
\[\lambda_i = \int_{a_i}^{b_i}\frac{2\sinh(\frac{t}{2})^2}{t}\, \df t.\]
\end{theorem}

We next consider the regime $\sqrt{g}\ll n(g) \ll g$. With this type of growth rate, the surfaces typically begin to possess many short closed geodesics. In fact, as can easily be seen in our proofs, if we do not shift the length window $[a,b]$ to be sufficiently close to zero with respect to the parameters $n$ and $g$, then $N_{g,n,[a,b]}$ do not even possess finite moments as $g\to\infty$. This is a similar phenomena to the $g=0$ setting in \cite{HideThomas2024} where the window was required to be rescaled by $\frac{1}{\sqrt{n}}$. 

Here we find that the correct scaling is $\frac{\sqrt{g}}{n(g)}$ and we obtain a different mean to the previous result for the limiting Poisson distributions. Notice that at this scale, the geodesics are all simple and disjoint by the collar lemma for hyperbolic surfaces (at least for $g$ large), and we find that the leading contribution to the moments comes from geodesics that form the boundary of an embedded pair of pants with two cusps. 

\begin{theorem}
   \label{th: stat de poisson dans le cas n^2/h -> +infty}
Let $k\in \N$ and
\[a_1 < b_1 < a_2 < b_2 < \ldots < a_k < b_k.\] Assume that
\[\lim_{g\to \infty}\frac{n(g)^2}{g} = \infty \qquad \text{and} \qquad \lim_{g\to \infty}\frac{n(g)}{g} = 0,\] then
\begin{equation*}
   \tuple{N_{g,n(g), \left[\frac{\sqrt{g}}{n(g)}a_1,\frac{\sqrt{g}}{n(g)}b_1\right]}}{N_{g,n(g), \left[\frac{\sqrt{g}}{n(g)}a_k,\frac{\sqrt{g}}{n(g)}b_k\right]}} \xrightarrow[g\to \infty]{(d)} \tuple{N_{[a_1, b_1]}}{N_{[a_k,b_k]}},
\end{equation*}
where $N_{[a_i, b_i]}$ are pairwise independent random variables following a Poisson distribution with parameter 
\[\lambda_i = \frac{b_i^2 - a_i^2}{32\pi^2}.\]

\end{theorem}

In our next result, we investigate the transitional regime when $n(g)$ and $\sqrt{g}$ are asymptotically comparable. 

\begin{theorem}
    \label{thm:stat de poisson dans le cas n^2/g -> alpha}
    Let $k \in \N$ and 
    \[a_1 < b_1 < a_2 < b_2 < \ldots < a_k < b_k < 2\, \mathrm{arcsinh}(1).\]  
    Assume that 
    \[\lim_{g\to \infty}\frac{n(g)^2}{g} = \alpha\in (0,\infty),\]
    and define
    \[\alpha_i = \frac{\alpha}{4\pi^2}\left(\cosh\left(\frac{b_i}{2}\right)-\cosh\left(\frac{a_i}{2}\right)\right) \quad \text{and} \quad \lambda_i = \int_{a_i}^{b_i} \frac{2\sinh^2\left(\frac{t}{2}\right)}{t} \, \df t.
     \]
    then
    \begin{equation*}
   \tuple{N_{g,n(g), \left[a_1,b_1\right]}(X)}{N_{g,n(g), \left[a_k,b_k\right]}(X)} \xrightarrow[g\to \infty]{(d)} \tuple{N_{[a_1, b_1]}}{N_{[a_k,b_k]}},
    \end{equation*}
    where $N_{[a_i, b_i]}$ are pairwise independent random variables following a Poisson distribution with parameter $\alpha_i+\lambda_i$.
       
\end{theorem}

As one might expect, the mean of the limiting Poisson distribution is a combination of the means in the previous two settings and the contributing curves to the moments are both the non-separating ones and those that separate a pair of pants with two cusps. Indeed, if we rescale the length window $[a,b]$ by $\frac{1}{\sqrt{\alpha}}$ (as done in Theorem \ref{th: stat de poisson dans le cas n^2/h -> +infty}) and consider $\alpha\to\infty$, then 
    $$\alpha_i = \frac{\alpha}{4\pi^2}\left(\cosh\left(\frac{b_i}{2\sqrt{\alpha}}\right)-\cosh\left(\frac{a_i}{2\sqrt{\alpha}}\right)\right) = \frac{b_i^2-a_i^2}{32\pi^2}+O\left(\frac{1}{\alpha}\right).$$
We also remark the restriction of the lengths to $<2\mathrm{arcsinh(1)}$. The reason for this is purely technical as it allows us to only consider the counting of simple and disjoint closed geodesics. We fully expect that to obtain an asymptotic Poisson distribution, we do not need to consider this restriction and we will investigate this in forthcoming work. Towards this, we readily obtain that the asymptotic joint factorial moments are finite without this length restriction (see Corollary \ref{cor: finite-moments} for a precise statement).

\begin{remark}
    Theorems \ref{th: stat de poisson dans le cas n^2/g -> 0} and \ref{th: stat de poisson dans le cas n^2/h -> +infty} can also be reframed in the language of Poisson point processes. For $X\in\mathcal{M}_{g,n}$, let $\xi_X^{(g)}$ denote the point measure associated with the multiset $\{\ell_\gamma(X)\}_\gamma$ (resp. $\{\frac{n}{\sqrt{g}}\ell_\gamma(X)\}_\gamma$) indexed over homotopy classes of closed curves on a surface of genus $g$ with $n$ cusps. Let $\xi^{(g)}$ denote the corresponding point process on $[0,\infty)$ with respect to $(\mathcal{M}_{g,n},\PP_{g,n})$. Theorem \ref{th: stat de poisson dans le cas n^2/g -> 0} (resp. Theorem \ref{th: stat de poisson dans le cas n^2/h -> +infty}) says that $\xi^{(g)}$ converges in distribution to a Poisson point process on $[0,\infty)$ with intensity $\frac{2\sinh^2\left(\frac{x}{2}\right)}{x}\df x$ (resp. $\frac{x}{16\pi^2}\df x$).
\end{remark}

\subsection{Proof ideas and further discussion}

As in \cite{MirzakhaniPetri2021} and \cite{HideThomas2024}, the proof of Theorems \ref{th: stat de poisson dans le cas n^2/g -> 0}, \ref{th: stat de poisson dans le cas n^2/h -> +infty} and \ref{thm:stat de poisson dans le cas n^2/g -> alpha} employs the method of factorial moments (see Theorem \ref{th : Method of factorial moments} below). This reduces the proof to counting the expected number of $k$-tuples of closed geodesics whose components are distinct and lengths live in one of the length windows.

When the components of the tuple are simple and do not intersect each other, their expected number is amenable for computation with Mirzakhani's integration formula (see Lemma \ref{lemma: Mirzakhani integration formula} below). This formula passes the counting to a computation of integrals involving Weil-Petersson volumes of moduli space of bordered hyperbolic surfaces. Many expansions of these volumes have been computed in the literature previously (see notably \cite{Mirzakhani2007-1,Mirzakhani2013,MirzakhaniZograf2024,WuXue2022,HideMaceraThomas2025,AnantharamanMonk2022,Hi.Th2025}). 

For our purposes, we require an extension of the results of Mirzakhani and Zograf \cite{MirzakhaniZograf2024} to our asymptotic $n^2(g)/g$ regimes (notably when $n^2(g)/g\to\infty$). To this end, we prove the following (the reader may refer to Section \ref{sec: background} for the notation used for intersection numbers and Weil-Petersson volumes if they are not familiar).

\begin{theorem}
    \label{th: developpement premier ordre comparaison classique}
    Let $n: \N \to \N$ be a function, satisfying
    \[\lim_{g \to \infty}\frac{n(g)}{g} = 0.\] 
    The following hold.
    \begin{enumerate}[(i)]
        \item Let $n_0\in\N$ be such that $n_0\leq n(g)$ for all sufficiently large $g$ and $\tuple{d_1}{d_{n_0}} \in \N^{n_0}$. Then,
            \begin{equation}
                \left(1- \frac{[\tau_{d_1}, \, \cdots,\, \tau_{d_{n_0}}, \tau_{0}^{n(g)-n_0}]_{g,n(g)}}{V_{g,n(g)}}\right) = \frac{\left(2|d|-n_0\right)n}{2\pi^2g} + O\!\left(|d|^4\left(\frac{n(g)^2}{g^2} + \frac{1}{g}\right)\right),
            \end{equation}
        \item
        \begin{equation}
            \label{eq: expansion 2 ème ordre V_g,n/V_g,n+1}
            \frac{4\pi^2(2g+n(g)-2)V_{g,n(g)}}{V_{g, n(g)+1}} = 1 + \frac{n(g)}{\pi^2g}+ O\!\left(\frac{n(g)^2}{g^2} + \frac{1}{g}\right),
        \end{equation}
        \item \begin{equation}
        \label{eq: expansion 3 ème ordre V_g,n/V_g,n+1}
            \frac{V_{g-1, n(g)+4}}{V_{g,n(g)+2}} = 1 - \frac{2n(g)}{\pi^2g} + O\!\left(\frac{n(g)^2}{g^2} + \frac{1}{g}\right).
        \end{equation}
    \end{enumerate}
\end{theorem} 
\begin{remark}
    Note that the first order terms in the previous theorem agree with those obtained by Mirzakhani and Zograf \cite{MirzakhaniZograf2024} for fixed $n$, the crucial extension is that we have explicit dependency on both $n$ and $g$ in the error term. 
    
    We could also develop the expansion further in powers of $n/g$. This was in fact carried out in \cite[Theorem 4.2]{HideMaceraThomas2025}, where a full asymptotic expansion was obtained; however the coefficients of the expansion were not made explicit (only bounded) as the spectral results obtained there did not require them. 
\end{remark}

Further to this, we also require a bound on products of Weil-Petersson volumes of unbordered moduli spaces summed over varying choices of genera and number of cusps. Many iterations of results of this form have also previously appeared in the literature (see for example \cite[Lemma 3.4]{Mirzakhani2013}, \cite[Lemma 6.3]{GLMST}, \cite[Lemma 22]{NieWuXue2022}, \cite[Lemma A.4]{Hide2023}) but again we require an extension that is valid up to $n(g)/g\to 0$. Our result is stated very general for suitability for application to a variety of problems in this asymptotic regime.

\begin{theorem}[Theorem \ref{th: somme d'un produit quelconque de V_gini} below]
Assume that
    $$\lim_{g\to\infty}\frac{n(g)}{g}=0.$$
   Let $q \geq 2$ be an integer, $(a_1,\, \cdots,\, a_q) \in \N^q$, $(s_1,\, \cdots,\, s_q) \in \left(\Z_{\geq 0}\setminus \{0,1,2\}\right)^q$ and $m \in \Z_{\geq 0}$. Then,
   \begin{equation}
      \dsum_{\substack{g_1 + \cdots + g_q = g \\ \bigsqcup I_i = \ens{1}{n}\\ 2g_i + \left|I_i\right| +a_i \geq s_i}}
      \frac{\prod_{i=1}^{q}V_{g_i, \, \left|I_i\right| +a_i}}{V_{g, \, n +m}} 
      = O_q\!\left(\dsum_{i = 1}^{q}\frac{1+n^{\sum_{\substack{j=1\\j \neq i}}^{q} s_j - a_j}}{g^{m - a_i + \sum_{\substack{j=1\\j \neq i}}^{q}s_j - a_j}}\right).
   \end{equation}
\end{theorem}

With these estimates at hand, we are able to readily compute the contributions of the simple and non-intersecting tuples to the factorial moments and identify which classes of tuples give rise to the leading orders.

The remaining task is to deal with tuples in which some of their components are either non-simple or have intersection between each other. Note that due to the restriction on the lengths that we consider and the collar lemma in hyperbolic geometry, this is only a concern for Theorem \ref{th: stat de poisson dans le cas n^2/g -> 0}. In this case, we can readily use very broad estimates (and simple extensions to the non-compact setting) that were present in the work of Mirzakhani and Petri (see \cite[Proposition 4.5]{MirzakhaniPetri2021}) to show that they only contribute subdominant terms to the factorial moments.

If we were to remove the $2\mathrm{arcsinh}(1)$ restriction from Theorem \ref{thm:stat de poisson dans le cas n^2/g -> alpha} we would also need to deal with the contribution from these tuples. In fact, the same arguments used for showing they are subdominant in the proof of Theorem \ref{th: stat de poisson dans le cas n^2/g -> 0} are sufficient to conclude that upon removing the $2\mathrm{arcsinh}(1)$ restriction, the factorial moments remain finite in the limit (see Corollary \ref{cor: finite-moments}). 

The difference in the $n(g)^2/g\to\alpha$ regime is that we do in fact expect non-simple and intersecting tuples to have a leading order contribution. Indeed, the core idea for bounding their contribution in Theorem \ref{th: stat de poisson dans le cas n^2/g -> 0} is to show that the intersections lead to the tuple of curves having a regular neighbourhood (retractable metric neighbourhood) that is a strict subsurface of the entire surface. Therefore, they give rise to a tuple of simple closed curves that together separate the surface which is a non-leading order contribution. 

However, when $n(g)^2/g\to\alpha$, we know that some of the separating curves do give leading order contributions, namely those that separate pants with two cusps from the surface. We therefore expect that tuples whose regular neighbourhoods are disjoint unions of these types of pants still give leading order contributions and we will elaborate on this in a future work.

\subsection{Corollaries on systoles and small eigenvalues}
The Poisson statistics results provide us some understanding on the size of the systole (the length of a shortest closed geodesic) and the number of small Laplacian eigenvalues on a hyperbolic surface.

For the regime $n(g)\ll\sqrt{g}$, we are able to deduce from using identical arguments to Mirzakhani and Petri (see \cite[Theorem 5.1]{MirzakhaniPetri2021}) that the expected size of the systole matches that of the closed surface setting as $g\to\infty$.
    \begin{corollary}
    \label{cor: systole}
    Assume that
\[\lim_{g\to \infty}\frac{n(g)^2}{g} = 0.\]
Then,
    $$\lim_{g\to\infty}\EE\left(\mathrm{systole}(X)\right)= \int_{0}^\infty e^{-\lambda_{[0,x]}}\df x=1.61498...,$$
where
    $$\lambda_{[0,x]}=\int_0^x\frac{2\sinh^2\left(\frac{t}{2}\right)}{t}\df t.$$
\end{corollary}
To obtain the expected systole size in the transitional regime $n(g)\asymp \sqrt{g}$ using the same arguments, we would need to remove the $2\mathrm{arcsinh}(1)$ length restriction. 

For the regime $\sqrt{g}\ll n(g)\ll g$, the scaling leads us to expect that the expected systole size will be of the order $\sqrt{g}/n(g)$. In fact, much of the same argument used to prove Corollary \ref{cor: systole} will follow through to show that
    $$\lim_{g\to\infty}\frac{n(g)\EE_{g,n}(\mathrm{systole}(X))}{\sqrt{g}}=\lim_{g\to\infty}\int_0^\infty \PP_{g,n}\left(N_{g,n(g),\left[0,\frac{\sqrt{g}}{n(g)}x\right]}(X)=0\right)\df x.$$
If we could justify the change of integral and limit, then by Theorem \ref{th: stat de poisson dans le cas n^2/h -> +infty} we would obtain
    $$\lim_{g\to\infty}\frac{n(g)\EE_{g,n}(\mathrm{systole}(X))}{\sqrt{g}}=\int_0^\infty \exp\left(-\frac{x^2}{32\pi^2}\right)\df x = 2\sqrt{2}\pi^{\frac{3}{2}}.$$
However, we were unable to obtain sufficient control over the random variable $N_{g,n(g),\left[0,\frac{\sqrt{g}}{n(g)}x\right]}$ to enable this. Instead, we can prove the following tail estimate.

\begin{proposition}
\label{prop: systole2}
    Assume that
        \[\lim_{g \to \infty}\frac{n(g)}{g} = 0 \qquad \text{and} \qquad \lim_{g \to \infty}n(g) = \infty.\]
    There exists a constant $C>0$ such that for any sequence $x_g$ satisfying 
        $$\lim_{g\to\infty}\frac{\sqrt{g}x_g}{n(g)}=0,$$
    we have
        $$\PP\left(\mathrm{systole}(X)>\frac{x_g\sqrt{g}}{n}\right)\leq C\left(\frac{1}{x_g^2}+x_g^2\left(\frac{1}{n}+\frac{g}{n^2}\right)+\frac{n}{g}\right),$$
    for $g$ sufficiently large.
\end{proposition}

In the regime $\sqrt{g}\ll n \ll g$, the existence of many very short closed geodesics that separate off two cusps as given by Theorem \ref{th: stat de poisson dans le cas n^2/h -> +infty}, results in the existence of many Laplacian eigenvalues close to zero. This is in contrast to what is known about the spectral gap in the other regimes that we consider (see Section \ref{sec: literature}).

Recall that for surfaces $X\in\mathcal{M}_{g,n}$ with $n\geq 1$, the Laplacian spectrum consists of an absolutely continuous part equal to $[\frac{1}{4},\infty)$, a simple eigenvalue at $0$ and then possibly discrete eigenvalues in $(0,\infty)$. The eigenvalues in $(0,\frac{1}{4})$ are called \textit{small eigenvalues} and their cardinality is bounded topologically by $2g+n-3$ by the work of Otal and Rosas \cite{Ot.Ro2009}. 

Using similar ideas as the proof of Theorem \ref{th: stat de poisson dans le cas n^2/h -> +infty} along with some standard min-max arguments (see \cite{HideThomas2024,Hi.Th2025}), we are able to show that in the $\sqrt{g}\ll n\ll g$ regime with probability tending to $1$ as $g\to\infty$, there is a lower bound on the number of small eigenvalues on a hyperbolic surface that is almost the same order as the Otal-Rosas upper bound.
\begin{theorem}
\label{thm: eigenavlues}
    Assume that
\[\lim_{g\to \infty}\frac{n(g)^2}{g} = \infty \qquad \text{and} \qquad \lim_{g\to \infty}\frac{n(g)}{g} = 0.\]   
    Then, for any $\varepsilon>0$, there exists a constant $C(\varepsilon)>0$ such that a Weil-Petersson random surface $X\in\mathcal{M}_{g,n}$ has at least $C(\varepsilon)\frac{n^2}{g}$ eigenvalues below $\varepsilon$ with probability tending to $1$ as $g\to\infty$.
\end{theorem}
\subsection{Comparison to related literature} 
\label{sec: literature}
As we noted, our work is inspired by the previous Poisson statistics results for lengths of closed geodesics obtained by Mirzakhani and Petri \cite{MirzakhaniPetri2021} and Hide and the second named author \cite{HideThomas2024} and the use of factorial moments in the proof originates in those works. The work of Mirzakhani and Petri has also inspired several other results in related settings. For example, Poisson statistics for counting closed curves on random planar and metric maps have appeared in the works \cite{BarazerGiacchettoLiu2023, JansonLouf2021, JansonLouf2021b} and for so-called tight closed geodesics on hyperbolic surfaces in \cite{BuddLions2025}. 

The disparity in the behaviour of the geodesics in the various regimes we consider is a direct consequence of the asymptotics of Weil-Petersson volumes of moduli spaces. When $n\ll \sqrt{g}$, these were considered in the seminal works of Mirzakhani \cite{Mirzakhani2007-1,Mirzakhani2007-2,Mirzakhani2013} and Mirzakhani and Zograf \cite{MirzakhaniZograf2024} with related results featuring in \cite{WuXue2022,LipnowskiWright2024,NieWuXue2022}. Much less has been understood in regimes outside of this. When $g$ is fixed and instead $n\to\infty$ asymptotics are considered, there is work by Zograf \cite{Zo1993} and Manin and Zograf \cite{Ma.Zo2000} and then later of Hide and the second named author \cite{Hi.Th2025}. When $n\ll g$ in general, there is the recent work \cite{HideMaceraThomas2025} of which our Theorems \ref{th: developpement premier ordre comparaison classique} and \ref{th: somme d'un produit quelconque de V_gini} complement. 

A lot of the recent work on understanding the geodesics in these different regimes has been undertaken in the context of studying the spectral gap of the Laplacian and the existence of eigenvalues below $\frac{1}{4}$. For typical closed surfaces (where the spectrum is purely discrete and there are infinitely many eigenvalues), it was long conjectured that the spectral gap is close to $\frac{1}{4}$, the bottom of the spectrum of their universal cover, $\mathbb{H}$. The existence of a sequence of closed surfaces of growing genus with $\lambda_1\to\frac{1}{4}$ was first proven by Hide and Magee \cite{Hi.Ma2023} as part of a broader result proving that random $n$-covers of cusped surfaces share the same property (with respect to their new eigenvalues) with probability tending to $1$ as $n\to\infty$.

In the Weil-Petersson random model, there were a series of works \cite{WuXue2022,LipnowskiWright2024,AnantharamanMonk2022} leading to a near optimal spectral gap of $\lambda_1>\frac{1}{4}-\varepsilon$ by Anantharaman and Monk \cite{anantharaman2025friedman} with probability tending to $1$ as $g\to\infty$. The state of the art in this regard is due to Hide, Macera and the second named author \cite{HideMaceraThomas2025} where a variant of the polynomial method \cite{chen2026new} was used to obtain a new proof leading to an explicit polynomial rate in $g$ for $\varepsilon$. We refer to \cite{magee2022random,magee2025strong,hide2025spectral,shen2025nearly} for other random models.

As with the behaviour of the geodesics, the near optimal spectral gap is expected to persist up to $n\ll\sqrt{g}$. This was first investigated by Hide \cite{Hide2023} and further by Shen and Wu in \cite{Sh.Wu2025}. The current state of the art is due to He, Wu and Xue \cite{he2026uniform} with a gap of $\frac{5}{16}-\varepsilon$.

When $\sqrt{g}\ll n\ll g$ we see in Theorem \ref{thm: eigenavlues} that this breaks down in a major way, with the existence of many small eigenvalues with high probability. The existence of one small eigenvalue in this regime was observed by Shen and Wu in \cite{Sh.Wu2025}. Our result is inspired by the analogous phenomena occurring when $g$ is fixed and $n\to\infty$ in the works of Hide and the second named author \cite{HideThomas2024,Hi.Th2025}. Perhaps more striking is that the existence of many (linear in volume) eigenvalues persists when $n>c g$ deterministically \cite{hide2024small} (the existence of at least one had been known for some time due to Zograf \cite{zograf1987small}).

The transitional regime $n^2\asymp g$ is also interesting and has been investigated before by Shen and Wu \cite{Sh.Wu2025}. They prove that there are asymptotically positive probability (depending upon $\alpha$) bounds for the Cheeger constant.

\subsection*{Acknowledgements}
This work was carried out during M.G.'s time as a research intern at Durham University from August 2025 to January 2026 and, along with future work on the $n^2\asymp g$ setting, will form M.G.'s master's dissertation. J.T. has received funding from the Leverhulme Trust through a Leverhulme Early Career Fellowship (Grant No. ECF-2024-440).

\subsection*{AI Disclosure}
AI was not used at any stage of this work.

\section{Definitions and background}
\label{sec: background}
In this section, we recall the main background notions for the Weil-Petersson random hyperbolic surface model as well as some preliminaries on Poisson random variables. A more in depth overview for the Teichm\"{u}ller theory can be found in \cite{Wright2019MirzakhaniSurvey}.

\subsection{Moduli space}
Let $g,n$ be non-negative integers such that $2g+n-2>0$. Let $\Sigma_{g,n}$ be a surface with genus $g$ and $n$ boundary components. A \textit{marked surface} of signature $(g,n)$ is a pair $(X, \phi)$, where $X$ is a hyperbolic surface with genus $g$ and $n$ geodesic boundary components, and $\phi:\Sigma_{g,n}\to X$ is an orientation-preserving diffeomorphism.

Let $L = \tuple{L_1}{L_n}\in \R_+^n$, we denote the Teichm\"{u}ller space of genus $g$ and $n$ labelled boundaries with lengths given by $L$ by
\[
\mathcal{T}_{g,n}(L_1,\cdots,L_n) \;\overset{\text{def}}{=}\;
\left\{
        \begin{array}{c}
        \text{Marked  surfaces } (X,\phi)\ \text{of signature } (g,n), \\[4pt]
        \text{ with labelled }\text{totally geodesic boundary components} \\[4pt]
        (\beta_1,\dots,\beta_n)\ \text{of lengths } (L_1,\dots,L_n)
        \end{array}
\right\} \big/ \sim,
\]
where $\sim$ is the equivalence relation identifying marked surfaces $(X_1, \phi_1)$ and $(X_2, \phi_2)$ if there exists an isometry $\psi:X_1\to X_2$ such that $\phi_{2}^{-1} \circ \psi \circ \phi_{1} $ is isotopic to the identity on $\Sigma_{g,n}$.

Now, let $\text{Homeo}^+(\Sigma_{g,n})$ denote the group of orientation-preserving homeomorphisms on $\Sigma_{g,n}$ which fix the boundary components set-wise. Let $\text{Homeo}_{0}^+(\Sigma_{g,n})$ be the subgroup of $\text{Homeo}^+(\Sigma_{g,n})$ of maps that are isotopic to identity. The mapping class group is then defined by
\[\text{MCG}_{g,n} \overset{\text{def}}{=} \text{Homeo}^+(\Sigma_{g,n})/\text{Homeo}_{0}^+(\Sigma_{g,n}).\]
The mapping class group acts on Teichmüller space through changing the marking
\[[\psi] \cdot [X,\phi] = \left[X, \phi \circ \psi^{-1}\right],\]
where $[X, \phi] \in \mathcal{T}_{g,n}(L)$ and $[\psi] \in \text{MCG}_{g,n}$.
The moduli space then arises as the quotient
\[\mathcal{M}_{g,n}(L) \overset{\text{def}}{=} \mathcal{T}_{g,n}(L)/ \text{MCG}_{g,n}.\]\\
The Teichmüller space $\mathcal{T}_{g,n}(L)$ carries a symplectic structure known as the Weil-Petersson symplectic form (see \cite{GOLDMAN1984200}) which we denote by $\omega_{\text{WP}}$. This form is invariant under the action of the mapping class group and hence descends to $\mathcal{M}_{g,n}(L)$ along with the volume form
\[\mathrm{d}\text{Vol}_{\text{WP},g,n} \overset{\text{def}}{=} \frac{1}{(3g+n-3)!} \bigwedge_{k=1}^{3g+n-3} \omega_{\text{WP}}.\]
On $\mathcal{M}_{g,n}(L)$, the volume form is finite (see \cite{Wolpert85}) and we denote its volume by $V_{g,n}(L)$. If \[L = \tuple{0}{0} = \mathbf{0},\] we denote $V_{g,n}(\mathbf{0})$ by $V_{g,n}$. If $k\geq n$ and $x_1,\ldots,x_k>0$ then we will also write $V_{g,n}(x_1,\ldots,x_k)$ for $V_{g,n}(x_1,\ldots,x_k,0,\ldots,0)$ where there are $n-k$ zeroes. We define the Weil-Petersson probability measure by normalizing the volume:  
\[\PP_{g,n}(\mathcal{A}) = \frac{1}{V_{g,n}(L)} \int_{\mathcal{M}_{g,n}(L)}\mathds{1}_{\mathcal{A}}(X)\, \mathrm{d}X,\]
where $\mathcal{A}$ is a measurable subset of $\mathcal{M}_{g,n}$ and we write $\mathrm{d}X$ as shorthand for $\mathrm{d}\text{Vol}_{\text{WP},g,n}(X)$. We will write $\EE_{g,n}$ for the expectation with respect to $\PP_{g,n}$. 

\subsection{Mirzakhani's integration formula}
In \cite{Mirzakhani2007-2}, Mirzakhani developed an integration formula for certain types of functions with respect to the Weil-Petersson volume form on the moduli space. The functions involve the lengths of multicurves on the surfaces. Recall the following definition. 
\begin{definition}
We define a $k$-multicurve as a $k$-tuple, $\Gamma = \tuple{\gamma_1}{\gamma_k}$, of non-homotopic non-peripheral, simple, closed, oriented curves on $\Sigma_{g,n}$ that are pairwise disjoint along with any curves in their free homotopy class. We denote by $[\Gamma]$ the free homotopy class of $\Gamma$.
\end{definition} 

Let $k$ be a positive integer, $[\psi]$ a mapping class group  element and $\Gamma= \tuple{\gamma_1}{\gamma_k}$ a $k$-multicurve. The mapping class group acts on $[\Gamma]$ by
\[[\psi]\cdot [\Gamma] = \left[\tuple{\psi\circ\gamma_1}{\psi\circ\gamma_k} \right]. \] 
Let
\[\mathcal{O}_{\Gamma} = \{[\psi]\cdot [\Gamma]\ |\ \psi \in \text{MCG}_{g,n} \},\]
denote the orbit of $\Gamma$ under this mapping class group action.
Let $\gamma$ be a closed curve on $\Sigma_{g,n}$ we can define the map $l_{\gamma}$ on the Teichm\"{u}ller space by associating to a pair $[X,f]$, the length of the unique closed geodesic in the free homotopy class of $f(\gamma)$ on $X$. Now let $F:\R^k \to\R$, we define the function $F_{\Gamma}$ by 
\[\begin{array}{rclc}
F_{\Gamma}: & \mathcal{M}_{g,n} & \longrightarrow & \mathbb{R^+} \\
& X & \longmapsto & \dsum_{\tuple{\alpha_1}{\alpha_k} \in \mathcal{O}_{\Gamma}} f\tuple{l_{\alpha_1}(X)}{l_{\alpha_k}(X)}
\end{array}.\]
where we use any fixed marking on $X$. This is well-defined on the moduli space due to the summation over the entire mapping class group orbit of $\Gamma$. 

Now, let $\Sigma_{g,n}(\Gamma)$ denote the (possibly disconnected) surface after cutting $\Sigma_{g,n}$ along $\Gamma$. That is,
    \[\Sigma_{g,n}(\Gamma) = \bigsqcup_{i=1}^s \Sigma_{g_i, n_i +b_i}.\]
for some $s \in \N$ and $\tuple{g_1,\,n_1,\, b_1}{g_s,\,n_s, \,b_s}$ non-negative integers where the $n_i$ denote the number of cusps from the original surface and $b_i$ the new boundaries formed from cutting along the curves in $\Gamma$.

Given $x = \tuple{x_1}{x_k}\in\mathbb{R}^k_{\geq 0}$, we assign the length $x_i$ to the curve $\gamma_i$ in $\Gamma$. Let $x^{(i)}$ be the tuple of coordinates $x_j$ corresponding to the lengths of the boundaries of the component $\Sigma_{g_i,n_i+b_i}$ formed by the curves in $\Gamma$ that make up its boundary. Thus each $x_j$ appears exactly twice among the tuples $x^{(i)}$, once for each side of the curve $\gamma_j$. We denote the volume of the moduli space corresponding to $\Sigma_{g,n}(\Gamma)$ by 
\[V_{g,n}(\Gamma, x) = \prod_{i=1}^{s} V_{g_i, n_i+b_i}\left(x^{(i)}\right).\]

Finally, denote by $M(\Gamma)$ to be the number of one-handles (genus $1$ and $1$ boundary component) that $\Gamma$ separates off. Mirzakhani's integration formula is then stated as follows.
\begin{lemma}[\cite{Mirzakhani2007-2}, Lemma 7.3]
\label{lemma: Mirzakhani integration formula}
For any function $F : \mathbb{R}^k \to \mathbb{R}_{+}$ and multicurve $\Gamma = (\gamma_1, \dots, \gamma_k)$,
the expectation of $F_{\Gamma}$ over $\mathcal{M}_{g,n}$ with respect to the Weil-Petersson volume form is given by
\[
\int_{\mathcal{M}_{g,n}} F_{\Gamma}(X) \, \mathrm{d}X
= 2^{-M(\Gamma)} \int_{\mathbb{R}_{+}^k} F(x_1,\ldots,x_k) \,V_{g,n}(\Gamma,x) \, x_1\cdots x_k \, \mathrm{d}x_1\cdots\mathrm{d}x_k.
\]
\end{lemma}

\subsection{Factorial moments}
\label{subsection : Factorial moments and particular notations}
In proving our main theorems, we will employ the method of factorial moments which we recall are defined as follows.
\begin{definition}
    Let $(\Omega, \PP)$ be a probability space and $X$ a random variable on $(\Omega, \PP)$ with values in $\Z$. Let $k$ be a positive integer. The $k$-th falling factorial of $X$ is given by
    \[(X)_k = X(X-1)\cdots(X-k+1).\]
    The $k$-th factorial moment of $X$ is then
    $\EE\left((X)_k\right).$
   
\end{definition} 

Recall that a random variable $X$ follows a Poisson distribution with parameter $\lambda$, denoted $X\sim\mathcal{P}(\lambda)$ if for every non-negative integer $k$, 
     \[\PP(X=k) = e^{-\lambda}\frac{\lambda^k}{k!}.\]

The factorial moments of a Poisson distributed random variable $X\sim\mathcal{P}(\lambda)$ take the particularly simple form of $\EE\left((X)_k\right)=\lambda^k$. The method of factorial moments states that these moments characterize Poisson distributed random variables.

\begin{theorem}[Method of factorial moments {\cite[Theorem 1.2]{Bollobas_2001}}]
\label{th : Method of factorial moments}
Let $l \in \mathbb{N}$ and suppose that $(\Omega_n, \mathbb{P}_n)_{n \geq 1}$ are a 
sequence of probability spaces, and $X_{i,n} : \Omega_n \to \mathbb{N}$ are a sequence of 
random variables for $i = 1, \dots, l$. For $k \in \mathbb{N}$,
Suppose that there exist $\lambda_i \in (0,\infty)$ for $i = 1,\dots, l$ such that for 
every $k_1, \dots, k_l \in \mathbb{N}$,
\[
\lim_{n \to \infty} \mathbb{E}_n \big[ (X_{1,n})_{k_1} \cdots (X_{l,n})_{k_l} \big]
= \lambda_1^{k_1} \cdots \lambda_l^{k_l}.
\]
Then, the random vector $(X_{1,n}, \dots, X_{l,n})$ converges in distribution to a 
random vector of independent Poisson distributed random variables with parameters $\lambda_i$.
\end{theorem}
In the rest of the paper, we will use the following notation adopted from Mirzakhani-Petri \cite{MirzakhaniPetri2021}. Given $a,b\in\mathbb{R}_{\geq 0}$ with $a<b$, define the random variable $N_{g,n, [a,b]}$ on $\mathcal{M}_{g,n}$ as 
    $$N_{g,n, [a,b]}(X) = \#\{\text{primitive, closed geodesics on $X$ with lengths in $[a,b]$}\}.$$
 
Given $r_1,\,\cdots,\, r_k\in\N$ and 
    $$0\leq a_1<b_1< a_2< \cdots <b_{k-1}< a_k <b_k,$$ 
the product of falling factorials, 
\[(N_{g,n, [a_1,b_1]})_{r_1} \cdots (N_{g,n, [a_k,b_k]})_{r_k},\]
counts the number of $\sum_ir_i$-tuples of primitive closed geodesics on $X$ where the first $r_1$ entries are distinct with lengths in $[a_1,b_1]$, the next $r_2$ entries are distinct with lengths in $[a_2,b_2]$ and so on. Note that the curves in the $r_i$ and $r_j$ tuples are also distinct because the length intervals do not overlap. 

We separate the product into a sum of two random variables, $Y_{g,n, r_1, \cdots, r_k}$ and $Y^\prime_{g,n, r_1, \cdots, r_k}$. The first counts the number of tuples where all of the entries are simple and pairwise disjoint and the second counts everything else, that is, the tuples of curves where some of the entries either intersect themselves or each other. 

\subsection{Volumes of moduli space}
Let $\overline{\mathcal{M}}_{g,n}$ denote the Deligne-Mumford compactification of $\mathcal{M}_{g,n}$. When $n>0$, there are $n$ line bundles $\mathcal{L}_i$ over $\overline{\mathcal{M}}_{g,n}$. For each $i$, $\mathcal{L}_i$ is the cotangent space at the $i$th marked point on $X \in \mathcal{M}_{g,n}$. Let $\psi_i$ denote the first Chern class of the $i$th line bundle over $\overline{\mathcal{M}}_{g,n}$.
For $d = \tuple{d_1}{d_n} \in \Z_{\geq 0}^n$ with
\[|d| = \sum_{i=1}^n d_i \leq 3g+n-3,\]
we define the intersection number
\[\left[\tau_{d_1}\cdots\tau_{d_n}\right]_{g,n} = \frac{2^{2|d|}\prod_{i=1}^{n}(2d_i +1)!!}{(3g +n-3- |d|)!}\int_{\overline{\mathcal{M}}_{g,n}}\psi^{d_1}\cdots\psi^{d_n}\omega_{\text{WP}}^{3g+n-3 -|d|}.\]  
If $|d|>3g+n-3$ then  $\left[\tau_{d_1}\cdots\tau_{d_n}\right]_{g,n}=0$ (see \cite{Mirzakhani2007-1} for more details). 

A remarkable result of Mirzakhani relates the volumes of moduli spaces of bordered surfaces to a polynomial whose coefficients are given by the intersection numbers. 
\begin{theorem}[{\cite[Theorem 1.1]{Mirzakhani2007-1}}]
    \label{Th: theorem 1.1 M_2007_1}
    Let $g,n \geq 0$ be such that $2g+n-2 >0$. Let $l_1,\, \ldots, \, l_n \in \R^{+}$ and $L=(l_1,\ldots,l_n)$. Then, 
    \[V_{g,n}(2L) = \sum_{\substack{d_1, \, \cdots\,d_n\\
                                    |d| \leq 3g+n-3}}
                                    \left[\tau_{d_1}\cdots\tau_{d_n}\right]_{g,n}
                \prod_{i=1}^n\frac{l_i^{2d_i}}{(2d_i +1)!}.                           
    \]        
\end{theorem}
    Thus the volumes of moduli spaces can be understood through the intersection numbers. These intersection numbers satisfy a number of recursive relations due to Liu-Xu \cite{LiuXU2009} and Mirzakhani \cite{Mirzakhani2007-2} that we will use later in this article. 
    \begin{lemma}
    \label{lemma: recursion_rels}
        Let $d= \tuple{d_1}{d_n} \in \Z_{\geq 0}^d$ such that
        \[|d| \leq 3g+n-3,\]
        and $d_0 = 3g+n-3-|d|$. The following recursive relations hold.
        \begin{flalign*}
            \left[\tau_0 \tau_1 \prod_{i=1}^n \tau_{d_i}\right]_{g,n+2} &=
            \left[\tau_0^4 \prod_{i=1}^n \tau_{d_i}\right]_{g-1,n+4}  \tag{I} \label{eq: rec I}\\
            &\quad + 6 \sum_{\substack{g_1+g_2=g \\ I \bigsqcup J = \{1,\dots,n\}}}
            \left[\tau_0^2 \prod_{i \in I} \tau_{d_i}\right]_{g_1,|I|+2} \cdot
            \left[\tau_0^2 \prod_{i \in J} \tau_{d_i}\right]_{g_2,|J|+2},
        \end{flalign*}
        \begin{flalign*}
            (2g - 2 + n)\left[\prod_{i=1}^n \tau_{d_i}\right]_{g,n} &=
            \frac{1}{2} \sum_{l=0}^{3g-3+n} (-1)^{l}(l+1) \frac{\pi^{2l}}{(2l+3)!}
            \left[\tau_{l+1} \prod_{i=1}^n \tau_{d_i}\right]_{g,n+1}. \tag{II} \label{eq: rec II}
        \end{flalign*}
        Moreover, if we set $a_0 =1/2$, and for $n \geq 1$,
        \[a_n = \zeta(2n)\left(1- 2^{1-2n} \right),\]
        where $\zeta$ is the Riemann zeta function then,
        \begin{flalign*}
            \left[\tau_{d_1} \cdots \tau_{d_n}\right]_{g,n} &=
            A_d + B_d + C_d,  \tag{III} \label{eq: rec III}
        \end{flalign*}
        where,
        \begin{align}
            A_d &= 8 \sum_{j=2}^n \sum_{l=0}^{d_0} (2d_j+1) a_l \,
            \left[\tau_{d_1 + d_j + l - 1} \prod_{i \neq 1,j} \tau_{d_i}\right]_{g,n-1} \label{eq: rec III A}\\
            B_d &= 16 \sum_{l=0}^{d_0} \sum_{k_1+k_2 = l + d_1 - 2} a_l \,
            \left[\tau_{k_1} \tau_{k_2} \prod_{i=2}^n \tau_{d_i}\right]_{g-1,n+1} \label{eq: rec III B}\\
            C_d &= 16 \sum_{\substack{g_1+g_2=g \\ I \bigsqcup J = \{2,\dots,n\}}}
            \sum_{l=0}^{d_0} \sum_{k_1+k_2 = l + d_1 - 2} a_l \,
            \left[\tau_{k_1} \prod_{i \in I} \tau_{d_i}\right]_{g_1,|I|+1} \cdot
            \left[\tau_{k_2} \prod_{i \in J} \tau_{d_i}\right]_{g_2,|J|+1}.\label{eq: rec III C}
        \end{align}
    \end{lemma}
    For further references, see \cite{Mirzakhani2013}. The coefficients $(a_n)_{n \in \N}$ satisfy the following estimates \cite{MirzakhaniZograf2024}.
    \begin{lemma}[{\cite[Lemma 2.1]{MirzakhaniZograf2024}}]
   \label{lemme:  aproxiamtion des coefficient a_i}
The sequence $\{ a_i \}_{i \geq 0}$ is increasing. Moreover,
\begin{enumerate}[i)]
    \item $\displaystyle \lim_{i \to \infty} a_i = 1$, and
    \[
    \sum_{i=0}^{\infty} (a_{i+1} - a_i) = \frac{1}{2}.
    \]
    \item $a_{i+1} - a_i$ has order $1 / 2^{2i}$, i.e., there exists $C > 0$ such that
    \[
    \frac{1}{C \cdot 2^{2i}} < a_{i+1} - a_i < \frac{C}{2^{2i}}.
    \]
    \item \[\displaystyle \sum_{i=0}^{\infty} i (a_{i+1} - a_i) = \frac{1}{4}.\]
   
\end{enumerate}
\end{lemma}
Using \eqref{eq: rec I} and \eqref{eq: rec II} Mirzakhani \cite{Mirzakhani2013} proved the following basic relations between volumes of moduli spaces.
\begin{lemma}
\label{Lemma 3.2 Mirzakhani2013}
For any $g,n \ge 0$ with $2g - 2 + n > 0$, the following estimates hold:
\[V_{g-1, n+4} \leq V_{g,n+2}\]
and
\begin{align}
\label{eq: bounded of Vg,n/Vg,n+1}
b_0 < (2g - 2 + n)\,\frac{V_{g,n}}{V_{g,n+1}} < b_1,
\end{align}
where $b_0$ and $b_1$ are universal constants independent of $g$ and $n$.\\

\end{lemma}
A better approximation of \eqref{eq: bounded of Vg,n/Vg,n+1} holds when the number of cusps $n$ is fixed and the genus $g$ tends to infinity.
\begin{theorem}[{\cite[Theorem 1.4]{Mirzakhani2013}}]
   \label{Th: quotient de V_g,n,V_g,n+1 avec n fixe}
   Let $n \in \Z_{\geq 0}$
   \[
\frac{4\pi^2(2g - 2 + n)V_{g,n}}{V_{g,n+1}} = 1 + O_n\!\left(\frac{1}{g}\right),
\]
as $g\to\infty$.
\end{theorem}
Finally, Mirzakhani and Zograf prove in \cite{MirzakhaniZograf2024} an asymptotic expansion for the volumes when $n$ is fixed.
\begin{theorem}[{\cite[Theorem 1.2]{MirzakhaniZograf2024}}]
\label{Th: expansion form n fix}
Let $n \in \Z_{\geq 0}$, then
   \[
V_{g,n} = \frac{C}{\sqrt{g}}(2g - 3 + n)! (4\pi^2)^{2g - 3 + n} \left(1 + O_n\left(\frac{1}{g}\right)\right),
\]
as $g \to \infty$.
\end{theorem} 

\section{Bounds and expansions for moduli space volumes in the regime $n\ll g$}
\label{sec: Computation of moduli space volume}
In this section, we prove two important results regarding the Weil-Petersson volumes of moduli spaces when the number of cusps $n\ll g$. The first is a bound on products of Weil-Petersson volumes. This serves as a generalization of Mirzakhani and Petri \cite[Lemma 3.2]{MirzakhaniPetri2021} when $n=0$ (see also \cite[Lemma 3.3]{Mirzakhani2013} and \cite[Lemma B.2]{LipnowskiWright2024}) and Hide \cite[Lemma A.5]{Hide2023} when $n\ll \sqrt{g}$. 

The second result that we prove here is Theorem \ref{th: developpement premier ordre comparaison classique}. This provides a first order expansion of certain intersection number and moduli space volume ratios for when $n\ll g$ with explicit dependency on $n$ in the error terms. This provides a generalization to the case where $n$ is fixed due to Mirzakhani and Zograf \cite[Theorem 1.4]{MirzakhaniZograf2024} and Lipnowski and Wright when $n\ll\sqrt{g}$ \cite[Theorem A.1]{LipnowskiWright2024}. We note that the results in \cite{HideMaceraThomas2025} can also be used to obtain higher order expansions.  

\subsection{Approximations for products of moduli space volumes}

Let $n:\N\to\N$ be a function satisfying the estimate
\[\lim_{g\to \infty} \frac{n(g)}{g} = 0.\]
Throughout, we will often omit the dependency on $g$ in the notation of the function $n$. In this section, we will prove the following bound on the product of Weil-Petersson volumes.
\begin{theorem}
    \label{th: somme d'un produit quelconque de V_gini}
   Let $q \geq 2$ be an integer, $(a_1,\, \cdots,\, a_q) \in \N^q$, $(s_1,\, \cdots,\, s_q) \in \left(\Z_{\geq 0}\setminus \{0,1,2\}\right)^q$ and $m \in \Z_{\geq 0}$. Then,
   \begin{equation}
      \label{eq: gestion d'un produit quelconque}
      \dsum_{\substack{g_1 + \cdots + g_q = g \\ \bigsqcup I_i = \ens{1}{n}\\ 2g_i + \left|I_i\right| +a_i \geq s_i}}
      \frac{\prod_{i=1}^{q}V_{g_i, \, \left|I_i\right| +a_i}}{V_{g, \, n +m}} 
      = O_q\!\left(\dsum_{i = 1}^{q}\frac{1+n^{\sum_{\substack{j=1\\j \neq i}}^{q} s_j - a_j}}{g^{m - a_i + \sum_{\substack{j=1\\j \neq i}}^{q}s_j - a_j}}\right).
   \end{equation}
\end{theorem}

\begin{remark}
    The product of volumes arises from using the Mirzakhani integration formula with a separating multicurve on a genus $g$ surface with $n$ cusps. The subsets $I_i$ appearing in the Theorem represent the subcollection of the (labelled) cusps that appear on the resulting $i$th component and the $a_i$ represent the number of its geodesic boundary components. For each $i$, the number $2-s_i$ represents an upper bound on the Euler characteristic of the $i$th component which is $\leq -1$ since the subsurfaces are hyperbolic, resulting in the requirement for $s_i\geq 3$.
\end{remark}

The proof of Theorem \ref{th: somme d'un produit quelconque de V_gini} will be done via an induction argument in $q\geq 2$. The base case corresponds to the following Lemma.
\begin{lemma}
   \label{lemme initialisation de rec produit des volumes}
   Let $(a_1, a_2) \in \N^2, \, (s_1, s_2) \in \left(\Z_{\geq 0} \setminus \{0, 1, 2\}\right)^{2}$ and $m \in \Z_{\geq 0}$. Then, 
   \begin{equation*}
      \dsum_{\substack{I \bigsqcup J = \ens{1}{n}\\ g_1 + g_2 = g \\ 2g_1 +|I| +a_1 \geq s_1 \\ 2g_2 +|J| +a_2 \geq s_2 }}
         \frac{V_{g_1, |I| +a_1} V_{g_2, |J| +a_2}}{V_{g,n+m}}
      = O\left( \frac{n^{s_1-a_1}}{g^{m+s_1-a_1-a_2}}+ \frac{n^{s_2-a_2}}{g^{m+s_2-a_1-a_2}}\right).
   \end{equation*}
\end{lemma}
To prove this lemma, we will subdivide the regime for $g_1$ depending on $g$ and $n$. We begin with when $g_1$ is large.
\begin{lemma}
    \label{lemma: initialisation rec somme des produits part 1}
    Let $(a_1, a_2) \in \N^2, \, (s_1, s_2) \in \left(\Z_{\geq 0} \setminus \{0, 1, 2\}\right)^{2}, \, m \in \Z_{\geq 0}$. Then,
    \begin{equation}
            \label{eq:  initialisation rec sommme des produits part 1}
            \dsum_{n_1 = 0}^n\dsum_{\substack{g_1 =\left\lfloor \frac{n-n_1 +a_2}{2}\right\rfloor +1\\ 2g_1 + n_1 +a_1 \geq s_1 \\ 2(g-g_1) + n-n_1 +a_2 \geq s_2}}^{\left\lfloor \frac{g}{2} \right\rfloor }
            \binom{n}{n_1} \frac{V_{g_1, n_1 +a_1} V_{g-g_1, n-n_1 +a_2}}{V_{g,n+m}}
            =O\!\left(\frac{1}{g^{m+s_1 -a_1 -a_2}}\right).
    \end{equation}
\end{lemma}
\begin{proof}
    For a pair $n_1 \in \ens{0}{n}$, $g_1 \in \ens{\left\lfloor \frac{n-n_1 +a_2}{2}\right\rfloor}{\left\lfloor\frac{g}{2}\right\rfloor}$, let $p \in \Z_{\geq 0}$ be so that $g_1 = p + \left\lfloor \frac{n-n_1 +a_2}{2}\right\rfloor$. Moreover, to simplify notation, define $n_2 :=n-n_1$. We now decompose the ratio $V_{g-g_1, n_2+a_2}/V_{g,n+m}$ as follows
\begin{align*}
      &\frac{V_{g-g_1, n_2 +a_2}}{V_{g,n+m}}\\
      =& \underbrace{\prod_{r = 0}^{\left\lfloor \frac{n_2 +a_2}{2}\right\rfloor -2} \frac{V_{g-g_1 +r , n_2+a_2 -2r}}{V_{g-g_1 +r +1, n_2+a_2 -2r-2}}}_{\text{\textcircled{a}}}
         \underbrace{\prod_{q=-1}^{p-1}\prod_{r=0}^{1} \frac{V_{g-g_1 + \left\lfloor \frac{n_2 +a_2}{2}\right\rfloor + q, n_2+a_2 -2\left\lfloor \frac{n_2 +a_2}{2}\right\rfloor + r+2}}{V_{g-g_1 + \left\lfloor \frac{n_2 +a_2}{2}\right\rfloor + q, n_2+a_2 -2\left\lfloor \frac{n_2 +a_2}{2}\right\rfloor + r+3}}}_{\text{\textcircled{b}}}\\
      & \underbrace{\prod_{q=-1}^{p-1}\frac{V_{g-g_1 + \left\lfloor \frac{n_2 +a_2}{2}\right\rfloor + q, n_2+a_2 -2\left\lfloor \frac{n_2 +a_2}{2}\right\rfloor + 4}}{V_{g-g_1 + \left\lfloor \frac{n_2 +a_2}{2}\right\rfloor + q+1, n_2+a_2 -2\left\lfloor \frac{n_2 +a_2}{2}\right\rfloor+2}}}_{\text{\textcircled{c}}}
         \underbrace{\prod_{r=0}^{n+m-3-n_2-a_2 +2\left\lfloor \frac{n_2 +a_2}{2}\right\rfloor} \frac{V_{g, n_2+a_2 -2\left\lfloor \frac{n_2 +a_2}{2}\right\rfloor + r+2}}{V_{g, n_2+a_2 -2\left\lfloor \frac{n_2 +a_2}{2}\right\rfloor + r+3}}}_{\text{\textcircled{d}}}.
\end{align*}
By Lemma \ref{Lemma 3.2 Mirzakhani2013}, \textcircled{a} and \textcircled{c} are bounded by $1$. Next, since
\[n_2+a_2 -2\left\lfloor \frac{n_2 +a_2}{2}\right\rfloor \in \{0,1\} \qquad \text{and} \qquad g-g_1 \geq \frac{g}{2},\]
we can use Theorem \ref{Th: quotient de V_g,n,V_g,n+1 avec n fixe} to find an $\alpha>1$ such that for any integer $q \geq -1$ and any $r \in \{0,1\}$ we have
\[\frac{V_{g-g_1 + \left\lfloor \frac{n_2 +a_2}{2}\right\rfloor + q, n_2+a_2 -2\left\lfloor \frac{n_2 +a_2}{2}\right\rfloor + r+2}}{V_{g-g_1 + \left\lfloor \frac{n_2 +a_2}{2}\right\rfloor + q, n_2+a_2 -2\left\lfloor \frac{n_2 +a_2}{2}\right\rfloor + r+3}}\leq \frac{\alpha}{4\pi^2 (2(g-g_1) +2q +n_2+a_2 +r  )}.\]
Moreover, we can choose $\alpha>1$ to be as close to 1 as long as we take $g$ sufficiently large. Thus,
\begin{align*}
      \text{\textcircled{b}}
      &\leq \prod_{q =-1 }^{p-1} \apar{\frac{\alpha}{4\pi^2}}^2 \frac{(2(g-g_1) +2q +n_2 +a_2 -1)!}{(2(g-g_1) +2q +n_2 +a_2 +1)!}\\
      &\leq  \apar{\frac{\alpha}{4\pi^2}}^{2p+2} \frac{(2(g-g_1) +n_2 +a_2 -3)!}{(2g +n_2 +a_2 - 2\left\lfloor \frac{n_2 +a_2}{2}\right\rfloor -1)!}.
   \end{align*}
By the second part of Lemma \ref{Lemma 3.2 Mirzakhani2013}, we bound \textcircled{d} by
   \[\text{\textcircled{d}} \leq b_{1}^{n+m-2-n_2-a_2 +2\left\lfloor \frac{n_2 +a_2}{2}\right\rfloor}\frac{(2g +n_2 +a_2 - 2\left\lfloor \frac{n_2 +a_2}{2}\right\rfloor -1)!}{(2g + n + m -3)!}.\]
Putting all of these bounds together, we obtain
   \begin{equation}
    \label{eq: 1 proof lemma initialisation}
    \frac{V_{g-g_1, n-n_1 +a_2}}{V_{g,n+m}} \leq \apar{\frac{\alpha}{4\pi^2}}^{2p+2}b_{1}^{n+m-2-n_2-a_2 +2\left\lfloor \frac{n_2 +a_2}{2}\right\rfloor} \frac{(2(g-g_1) +n_2 +a_2 -3)!}{(2g + n + m -3)!}.
   \end{equation}
By Lemma \ref{Lemma 3.2 Mirzakhani2013}, Theorem \ref{Th: expansion form n fix} and \eqref{eq: 1 proof lemma initialisation} we bound
\begin{align*}
   &\frac{V_{g_1, n_1 +a_1} V_{g-g_1, n-n_1 +a_2}}{V_{g,n+m}}\binom{n}{n_1}\\
   \leq& \binom{n}{n_1}\frac{V_{g_1 + \left\lfloor \frac{n_1 +a_1}{2}\right\rfloor-1, n_1 +a_1 -2 \left\lfloor \frac{n_1 +a_1}{2}\right\rfloor +2} V_{g-g_1, n-n_1 +a_2}}{V_{g,n+m}}\\
   \leq& C \binom{n}{n_1}(4\pi^2)^{2g_1+n_1+a_1 -3}(2g_1+n_1+a_1 -3)! b_{1}^{n+m}\left(\frac{\alpha}{4\pi^2}\right)^{2p }\frac{(2(g-g_1) +n_2 +a_2 -3)!}{(2g + n + m -3)!}\\
   \leq& C \binom{n}{n_1} (4\pi^2)^{n}b_{1}^{n}\frac{\alpha^{2p}}{g^{m+3-a_1-a_2}}\frac{1}{\binom{2g+n +a_1 +a_2-6}{2g_1 +n_1 +a_1 -3}},
\end{align*}
where the constant $C>0$ changes from line to line but is independent of $n$ and $g$. Now note that
\[\binom{n}{n_1} \leq 2^n,\] 
\[\apar{\frac{2g+n+a_1 +a_2 -6}{2g_1 + n_1 +a_1 -3}}^{2g_1 + n_1 +a_1 -3}\leq\binom{2g+n +a_1 +a_2-6}{2g_1 +n_1 +a_1 -3},\]
and $\alpha^{2p} \leq \alpha^{2g_1} \leq C\alpha^{2g_1 +n_1+ a_1 -3}$. Thus, 
\[\frac{V_{g_1, n_1 +a_1} V_{g-g_1, n-n_1 +a_2}}{V_{g,n+m}}\binom{n}{n_1} \leq C \frac{ (8\pi^2b_1)^{n}}{g^{m+3-a_1-a_2}}\apar{\frac{\alpha(2g_1 + n_1 +a_1 -3)}{2g+n+a_1 +a_2 -6}}^{2g_1 + n_1 +a_1 -3}.\]
To bound this, we use the following fact.\\
\\
\begin{fact}
\label{fact: convexity} Let $a\in\R_{>0}$ and let $f$ be the function defined by
\begin{equation}
    \label{eq: convexité de (x/a)^x}
    \begin{array}{rclc}
        f: & \R_{\geq 0} & \longrightarrow & \R_{\geq 0} \\
        & x & \longmapsto & \left(\frac{x}{a}\right)^x
    \end{array}.
\end{equation}
\end{fact}
\textit{Then, $f$ is convex and in particular, when restricting $f$ to a compact interval $I\subseteq \R_{\geq 0}$, the maximum is at an endpoint of $I$.}\\
\\
In our case, 
\[0\leq n_1 \leq n, \qquad \lfloor \frac{n-n_1+a_2}{2}\rfloor +1\leq g_1 \leq \lfloor\frac{g}{2}\rfloor \qquad \text{and} \qquad 2g_1 +n_1 +a_1 \geq s_1.\]
So $2g_1 + n_1 +a_1 -3\geq s_1-3$ and $2g_1+n_1+a_1-3\geq n+a_1+a_2-3$ and thus we can restrict to the interval
\[I = \left[\max\left(n + a_2 + a_1 -3, \, s_1-3 \right), \,g+n+a_1-3 \right].\]
In the case $2g_1 +n_1 +a_1- 3 = \max\left(n + a_2 + a_1 -3, \, s_1-3 \right)$, we obtain
\[\left(\frac{\alpha(2g_1 +n_1 +a_1 -3)}{2g +n +a_1 +a_2 - 6}\right)^{2g_1 +n_1 +a_1 -3} \leq \left(\frac{\alpha\max\left(n + a_2 + a_1 -3, \, s_1-3 \right)}{2g +n +a_1 +a_2-6}\right)^{\max\left(n + a_2 + a_1 -3, \, s_1 -3 \right)}.\]

When $2g_1 +n_1 +a_1- 3 = g+ n +a_1 -3$,
because for $g\geq g_0$ for some fixed $g_0$ we can take $\alpha$ to be smaller than $3/2$, we obtain
\[\left(\alpha\frac{g+ n + a_1 -3}{2g+n+a_1+a_2-6}\right)^{g+n +a_1 -3} \leq C\left(\frac{\alpha}{2}\right)^{g+n}e^{n} \leq Ce^{n} \left(\frac{3}{4}\right)^{g+n}.\]
Summing these two and separating out the case where $2g_1+n_1 +a_1 = s_1$, we see that \eqref{eq: initialisation rec sommme des produits part 1} is bounded by
\begin{align*}
    &\dsum_{n_1 = 0}^n\dsum_{\substack{g_1 =\left\lfloor \frac{n-n_1 +a_2}{2}\right\rfloor +1\\ 2g_1 + n_1 +a_1 \geq s_1 \\ 2(g-g_1) + n-n_1 +a_2 \geq s_2}}^{\left\lfloor \frac{g}{2} \right\rfloor }
            \binom{n}{n_1} \frac{V_{g_1, n_1 +a_1} V_{g-g_1, n-n_1 +a_2}}{V_{g,n+m}}\\
    &\leq C
    \dsum_{n_1 = 0}^n\frac{ (8\pi^2b_1)^{n}}{g^{m+3-a_1-a_2}} \left(
        \left(\alpha\frac{s_1-3}{2g +n +a_1 +a_2-6}\right)^{s_1 -3 }\right.\\
    &\left.\phantom{\dsum_{n_1 = 0}^n}+\dsum_{\substack{g_1 =\left\lfloor \frac{n-n_1 +a_2}{2}\right\rfloor +1\\ 2g_1 + n_1 +a_1 -3 \geq \max\left(n + a_2 + a_1 -3, \, s_1-2 \right) \\ 2(g-g_1) + n-n_1 +a_2 \geq s_2}}^{\left\lfloor \frac{g}{2} \right\rfloor }
    \left(\alpha\frac{\max\left(n + a_2 + a_1 -3, \, s_1-2 \right)}{2g +n +a_1 +a_2-6}\right)^{\max\left(n + a_2 + a_1 -3, \, s_1 -2 \right)} \right.\\
    &\phantom{ \dsum_{n_1 = 0}^n\dsum_{\substack{g_1 =\left\lfloor \frac{n-n_1 +a_2}{2}\right\rfloor +1\\ 2g_1 + n_1 +a_1 \geq s_1 \\ 2(g-g_1) + n-n_1 +a_2 \geq s_2}}^{\left\lfloor \frac{g}{2} \right\rfloor } \frac{ (8\pi^2b_1)^{n}}{g^{m+3-a_1-a_2}}} \left.
   + \frac{1}{e^n} \left(\frac{3}{4}\right)^{g+n}\right)\\
    &\leq 2C(8\pi^2b_1)^{n}\frac{1+n}{g^{m+2-a_1-a_2}}\left(\alpha\frac{\max\left(n + a_2 + a_1 -3, \, s_1-2 \right)}{2g +n +a_1 +a_2-6}\right)^{\max\left(n + a_2 + a_1 -3, \, s_1 -2 \right)}.
\end{align*}
To conclude, we distinguish between the case when $n$ is bounded or not.
If $n$ is bounded, $(8\pi^2b_1)^{n}$ can be absorbed as a constant and so we have the bound
\[\frac{1+n}{g^{m+2-a_1-a_2}}\left(\alpha\frac{\max\left(n + a_2 + a_1 -3, \, s_1-2 \right)}{2g +n +a_1 +a_2-6}\right)^{\max\left(n + a_2 + a_1 -3, \, s_1 -2 \right)} \leq C\frac{1}{g^{m +s_1 -a_1 -a_2}},\]
since if $n + a_1+a_2-3 \geq s_1-2$ then $1/g^{n + a_1+a_2-3} \leq 1/g^{s_1-2}$.
If $n$ is unbounded then, for $g$ large enough, we have
\begin{align*}
(8\pi^2b_1)^{n}\frac{1+n}{g^{m+3-a_1-a_2}}\left(\alpha\frac{n + a_2 + a_1 -4}{2g+n+a_1+a_2-6}\right)^{n+a_2 +a_1 -4} 
&\leq \frac{1+n}{g^{m+3-a_1-a_2}}\left(C\frac{1+n}{g}\right)^{n+a_1+a_2 -4} \\
&\leq C \frac{1}{g^{m+s_1-a_1 -a_2}},
\end{align*}
where in the last inequality we have used that because $n\to\infty$ and $n/g\to 0$ as $g\to\infty $, we have for $g$ large enough that
    $$\left(\frac{C(1+n)}{g}\right)^n\leq \frac{1}{g^k},$$
for any fixed $k$. Combining all of the estimates concludes the proof.
\end{proof}
We now deal with when $g_1$ is small.
\begin{lemma}
    \label{lemma initialisation rec somme des produits part 2}
    Let $(a_1, a_2) \in \N^2, \, (s_1, s_2) \in \left(\Z_{\geq 0} \setminus \{0, 1, 2\}\right)^{2}, \, m \in \N$. Then,
    \begin{equation}
            \label{eq:  initialisation rec sommme des produits part 2}
            \dsum_{n_1 = 0}^n\dsum_{\substack{g_1 = 0\\ 2g_1 + n_1 +a_1 \geq s_1 \\ 2(g-g_1) + n-n_1 +a_2 \geq s_2}}^{\left\lfloor \frac{n-n_1 +a_2}{2} \right\rfloor }
            \binom{n}{n_1} \frac{V_{g_1, n_1 +a_1} V_{g-g_1, n-n_1 +a_2}}{V_{g,n+m}} 
            = O\!\left(\frac{1+n^{s_1-a_1}}{g^{m+s_1-a_1-a_2}}\right).
        \end{equation}
\end{lemma}

\begin{proof}
 
The proof of this bound is similar to Lemma \ref{lemma: initialisation rec somme des produits part 1}, the main difference is in the application of Fact \ref{fact: convexity}. 

First, remark that because $g_1 \leq \left\lfloor \frac{n-n_1 +a_2}{2}\right\rfloor$ the condition $2(g-g_1)+n-n_1+a_2 \geq s_2$ holds whenever $g\geq \frac{s_2}{2}$ and hence for $g$ sufficiently large. So, making the change of variable $t_1 = 2g_1 +n_1 +a_1$ in the inner summation, extending the resulting sum up to $n+a_1+a_2$ and then interchanging summations, we can bound 
\begin{align*}
   &\dsum_{n_1 = 0}^n\dsum_{\substack{g_1 = 0\\ 2g_1 + n_1 +a_1 \geq s_1 \\ 2(g-g_1) + n-n_1 +a_2 \geq s_2}}^{\left\lfloor \frac{n-n_1 +a_2}{2} \right\rfloor }
            \binom{n}{n_1} \frac{V_{g_1, n_1 +a_1} V_{g-g_1, n-n_1 +a_2}}{V_{g,n+m}}\\
   &=\sum_{n_1=0}^n\sum_{\substack{t_1=n_1+a_1 \\ g_1 = \frac{t_1-a_1-n_1}{2} \in \Z}}^{2\left\lfloor \frac{n-n_1 +a_2}{2} \right\rfloor+n_1+a_1} \binom{n}{n_1} \frac{V_{g_1, n_1 +a_1} V_{g-g_1, n-n_1 +a_2}}{V_{g,n+m}}\1_{\{t_1\geq s_1\}}(t_1,n_1)\\
   &\leq \sum_{n_1=0}^n\sum_{\substack{t_1=\max\{s_1,n_1+a_1\} \\g_1 = \frac{t_1-a_1-n_1}{2} \in \Z}}^{n+a_1+a_2} \binom{n}{n_1} \frac{V_{g_1, n_1 +a_1} V_{g-g_1, n-n_1 +a_2}}{V_{g,n+m}} \\
   &=\underbrace{\sum_{n_1=0}^{s_1-a_1}\sum_{\substack{t_1=s_1 \\ g_1 = \frac{t_1-a_1-n_1}{2} \in \Z}}^{n+a_1}}_{(a)}+\underbrace{\sum_{n_1=s_1-a_1+1}^n\sum_{\substack{t_1=n_1+a_1 \\ g_1=t_1-n_1-a_1\in 2\Z}}^{n+a_1}}_{(b)}+\underbrace{\sum_{n_1=0}^n\sum_{\substack{t_1=n+a_1+1 \\ g_1 = \frac{t_1-a_1-n_1}{2} \in \Z}}^{n+a_1+a_2}}_{(c)},    
\end{align*}
where the summands on the final line are the same as on the previous line.
As in the proof of Lemma \ref{lemma: initialisation rec somme des produits part 1}, we bound the contribution of $V_{g-g_1, n-n_1+a_2}$ using Lemma \ref{Lemma 3.2 Mirzakhani2013},
\begin{align*}
   \frac{V_{g-g_1, n-n_1+a_2}}{V_{g,n+m}} 
   &=\prod_{r=0}^{g_1-1} \frac{V_{g-g_1+r, n-n_1+a_2 -2r}}{V_{g-g_1 +r+1, n-n_1+a_2 -2r-2}}\prod_{r=0}^{n_1+m-1+2g_1-a_2} \frac{V_{g, n-n_1+a_2 -2g_1+ r}}{V_{g, n-n_1+a_2-2g_1 +r+1}}\\
   &\leq \prod_{r=0}^{n_1+m-1+2g_1-a_2} \frac{V_{g, n-n_1-2g_1+a_2 + r}}{V_{g, n-n_1-2g_1+a_2 +r+1}}\\
   &\leq b_{1}^{n_1+m+2g_1-a_2}\frac{(2(g-g_1) +n-n_1 +a_2 -3)!}{(2g +n +m-3)!},
\end{align*}
and so, as in the proof of Lemma \ref{lemma: initialisation rec somme des produits part 1}, we obtain
\begin{equation}
    \label{eq: 2 proof lemma initialisation}
   \frac{V_{g_1, n_1 +a_1} V_{g-g_1, n-n_1 +a_2}}{V_{g,n+m}}\binom{n}{n_1}
   \leq C\binom{n}{n_1}\frac{(4\pi^2b_1)^{t_1}}{g^{m+3-a_1 -a_2}}\apar{\frac{t_1 -3}{2g+n+a_1 +a_2 -6}}^{t_1 -3}.
\end{equation}
For $(c)$, we can then repeat the method used in the proof of Lemma \ref{lemma: initialisation rec somme des produits part 1} to show that
\begin{equation}
    \label{eq: proof initialisation (2'')}
    (c) \leq C\frac{1}{g^{m +s_1 -a_1 -a_2}}.
\end{equation}

The main difference between bounding (a) and (b) lies in the upper bound used for the binomial coefficient. In the first case, we will use the bound 
\begin{equation}
    \label{eq: binom < n^k}
\binom{n}{n_1}\leq n^{n_1},
\end{equation}
and in the second, 
\begin{equation}
    \label{eq: binom < (en/k)^k}
\binom{n}{n_1} \leq \left(\frac{en}{n_1}\right)^{n_1}.
\end{equation}
Additionally, in the case of (a), we separate the case where $t_1 = s_1$ since this corresponds to the main leading order term.
So, using \eqref{eq: 2 proof lemma initialisation}, we get
\begin{align*}
    (a)\leq&C\dsum_{n_1 = 0}^{s_1-a_1} 
         \binom{n}{n_1}\frac{(4\pi^2b_1)^{s_1}}{g^{m+3-a_1-a_2}}\apar{\frac{s_1  -3}{2g+n+a_1 +a_2 -6}}^{s_1 -3}\\
         &+C\dsum_{n_1 = 0}^{s_1-a_1} \sum_{t_1 = s_1+1}^{n+a_1}
         \binom{n}{n_1}\frac{(4\pi^2b_1)^{t_1}}{g^{m+3-a_1 -a_2}}\apar{\frac{t_1 -3}{2g+n+a_1 +a_2 -6}}^{t_1 -3}.
\end{align*}
The first term is dealt with trivially, and the second uses Fact \ref{fact: convexity}, noting that for $g$ sufficiently large, the maximum is achieved when $t_1=s_1+1$ and there are at most $C(1+n)$ terms in the inner summation for some constant $C>0$ depending only on $a_1$ and $s_1$. We hence bound (a) by

\begin{align}
    \label{eq: proof initialisation (e)}
(a)
&\leq C\left(\frac{1}{g^{m+s_1-a_1-a_2}} + \frac{n}{g^{m+1+s_1-a_1-a_2}}\right) \dsum_{n_1 = 0}^{s_1-a_1} n^{n_1}
\leq C\frac{1+n^{s_1-a_1}}{g^{m+s_1-a_1-a_2}}. 
\end{align}
Finally for (b), we bound the volume term using \eqref{eq: 2 proof lemma initialisation}, separate out the $t_1 = n_1 +a_1$ term in the inner summation and then use Fact \ref{fact: convexity} to bound the remaining terms by their value at $t_1= n_1+a_1+1$ resulting in 

 \begin{align*}
    (b)
    \leq&C\dsum_{n_1 = s_1-a_1 +1}^{n} 
         \frac{1}{g^{m-a_2}}\apar{\frac{4\pi^2b_1en}{2g+n+a_1 +a_2 -6}}^{n_1}n_{1}^{a_1 -3}\\
      &+C\dsum_{n_1 = s_1-a_1 +1}^{n}
         \frac{n}{g^{m+1 -a_2}}\apar{\frac{4\pi^2b_1en}{2g+n+a_1 +a_2 -6}}^{n_1}n_{1}^{a_1 -2}
 \end{align*} 
To conclude, we use a new fact which is proven using elementary calculus.\\
\\
\begin{fact} 
\label{fact: fact2} Given $a,b,c\in\mathbb{N}$, the function $h$ defined by
\begin{equation}
    \label{eq: convexité de (n/g)^x x^c}
    \begin{array}{rclc}
        h: & [a, \, \infty) & \longrightarrow & \R_{>0} \\
        & x & \longmapsto & \left(\frac{bn}{g}\right)^x x^c
    \end{array}
\end{equation}
is decreasing whenever $a>\frac{c}{\log(g)-\log(bn)}$ and so the maximum is attained at $a$.
\end{fact}
We apply Fact \ref{fact: fact2} to 
\[\apar{\frac{4\pi^2b_1en}{2g+n+a_1 +a_2 -6}}^{n_1}n_{1}^{a_1 -2} \qquad \text{and } \qquad \apar{\frac{4\pi^2b_1en}{2g+n+a_1 +a_2 -6}}^{n_1}n_{1}^{a_1 -3}.\]
So for $g$ large enough, the maximum is obtained at $n_1 =s_1-a_1 +1$. Thus, 
\begin{equation}
    \label{eq: proof initialisation (f)}
    (b)
    \leq C\frac{n^{s_1-a_1 +1}}{g^{m+ 1 +s_1 -a_1-a_2}} +C\frac{n^{s_1-a_1 +2}}{g^{m+ 2 +s_1 -a_1-a_2}}
    \leq C\frac{1+n^{s_1-a_1}}{g^{m +s_1 -a_1-a_2}}
\end{equation}
The result now follows from combining \eqref{eq: proof initialisation (2'')}, \eqref{eq: proof initialisation (e)} and \eqref{eq: proof initialisation (f)}.    
\end{proof}

The proof of Lemma \ref{lemme initialisation de rec produit des volumes} is now a combination of Lemmas \ref{lemma: initialisation rec somme des produits part 1} and \ref{lemma initialisation rec somme des produits part 2}.
\begin{proof}[Proof of Lemma \ref{lemme initialisation de rec produit des volumes}]
We rewrite the summation as
\begin{align*}
    \dsum_{\substack{I \bigsqcup J =
    \ens{1}{n}\\ g_1 + g_2 = g \\ 2g_1 +|I| +a_1 \geq s_1 \\ 2g_2 +|J| +a_2 \geq s_2 }}
         \frac{V_{g_1, |I| +a_1} V_{g_2, |J| +a_2}}{V_{g,n+m}}
         =& \dsum_{g_1 = 0}^{\left\lfloor \frac{g}{2} \right\rfloor } \dsum_{\substack{n_1 = 0\\ 2g_1 + n_1 +a_1 \geq s_1 \\ 2(g-g_1) + n-n_1 +a_2 \geq s_2}}^n 
         \binom{n}{n_1} \frac{V_{g_1, n_1 +a_1} V_{g-g_1, n-n_1 +a_2}}{V_{g,n+m}}\\
       & + \dsum_{g_2 = 0}^{\left\lfloor \frac{g}{2} \right\rfloor - \epsilon} \dsum_{\substack{n_2 = 0\\ 2g_2 + n_2 +a_2 \geq s_2 \\ 2(g-g_2) + n-n_2 +a_1 \geq s_1}}^n 
         \binom{n}{n_2} \frac{V_{g-g_2, n-n_2 +a_1} V_{g_2, n_2 +a_2}}{V_{g,n+m}},\\
\end{align*}
where $\epsilon$ is $0$ if $g$ is odd and $1$ if $g$ is even. Because we are interested only by an upper bound, we can consider $\epsilon$ to be $0$ resulting in identical terms on the right hand side after interchanging the roles of $s_1$ with $s_2$, $a_1$ with $a_2$, $n_1$ with $n_2$ and $g_1$ with $g_2$.

To bound such a summation, we split it depending on whether $g_1$ is smaller or larger than $\left\lfloor \frac{n-n_1+a_2}{2}\right\rfloor$. The resulting summations are then precisely the terms bounded in Lemmas \ref{lemma: initialisation rec somme des produits part 1} and \ref{lemma initialisation rec somme des produits part 2} and hence the result follows.

\end{proof}
We now prove Theorem \ref{th: somme d'un produit quelconque de V_gini} by an induction over $q$, the number of volumes in the product. 

\begin{proof}[Proof of Theorem \ref{th: somme d'un produit quelconque de V_gini}]
Let $q\geq 2$ be an integer for which the result holds for a product of $q$ volumes. The base case for $q=2$ is Lemma \ref{lemme initialisation de rec produit des volumes}.

For a product of $q+1$ volumes, since the genera sum to $g$, at least one of the $g_i$ must be at least $g/(q+1)$ and so we may bound 
\begin{equation}
    \label{eq: th somme des produits rec. Decomposition de la somme principal}
    \dsum_{\substack{g_1 + \cdots + g_{q+1} = g \\ \bigsqcup I_i = \ens{1}{n}\\ 2g_i + \left|I_i\right| +a_i \geq s_i}}
      \frac{\prod_{i=1}^{q+1}V_{g_i, \, \left|I_i\right| +a_i}}{V_{g, \, n +m}}
   \leq \dsum_{j =1}^{q+1} \dsum_{\substack{g_1 + \cdots + g_{q+1} = g \\ \bigsqcup I_i = \ens{1}{n}\\ 2g_i + \left|I_i\right| +a_i \geq s_i\\ g_{j} \geq \lfloor \frac{g}{q+1} \rfloor}}
            \frac{\prod_{i=1}^{q+1}V_{g_i, \, \left|I_i\right| +a_i}}{V_{g, \, n+m}}
\end{equation}
The proceeding bounds will be independent of the index $j$ so to be concrete we will assume $j=q+1$ in the following. We will apply Lemma \ref{lemme initialisation de rec produit des volumes} to the summation  only over the indices involving $q$ and $q+1$, that is keeping the variables for the other indices fixed. We have
\begin{equation}
\label{eq: th sommm des produits rec. séparation des q-1 premiers termes et des deux autres}
\begin{aligned}         
   &\dsum_{\substack{g_1 + \cdots + g_{q+1} = g \\ \bigsqcup I_i = \ens{1}{n}\\ 2g_i + \left|I_i\right| +a_i \geq s_i\\ g_{q+1} \geq \lfloor \frac{g}{q+1} \rfloor}}
            \frac{\prod_{i=1}^{q+1}V_{g_i, \, \left|I_i\right| +a_i}}{V_{g, \, n+m}}\\
   \leq &\dsum_{\substack{n_1 + \cdots + n_{q-1} \leq n \\ g_1 + \cdots + g_{q-1} \leq g - \lfloor \frac{g}{q+1} \rfloor \\ 2g_i + n_i +a_i \geq s_i}}
            \frac{n!}{n_1! \cdots n_{q-1}! \left(n- \dsum_{i=1}^{q-1} n_i\right)!}\frac{\prod_{i=1}^{q-1}V_{g_i, \, n_i +a_i}}{V_{g, \, n+m}}\\
         &\cdot   \dsum_{n_{q} +n_{q+1} = n - \sum_{i =1}^{q-1} n_i} \dsum_{\substack{g_q + g_{q+1} = g - \sum_{i= 1}^{q-1} g_i\\ 2g_q + n_q +a_q \geq s_q\\ 2g_{q+1} +n_{q+1} +a_{q+1} \geq 2\lfloor\frac{g}{q+1} \rfloor }}
            \binom{n- \sum_{i=1}^{q-1} n_i}{n_{q+1}} V_{g_q, \, n_q +a_q} V_{g_{q+1}, \, n_{q+1} + a_{q+1}}.
   \end{aligned}
\end{equation}
Thus, by Lemma \ref{lemme initialisation de rec produit des volumes} applied with $n$ replaced by $n-\sum_{i=1}^{q-1} n_i$.
\begin{align*}
    &\dsum_{n_{q} +n_{q+1} = n - \sum_{i =1}^{q-1} n_i} \dsum_{\substack{g_q + g_{q+1} = g - \sum_{i= 1}^{q-1} g_i\\ 2g_q + n_q +a_q \geq s_q\\ 2g_{q+1} +n_{q+1} +a_{q+1} \geq 2\lfloor\frac{g}{q+1} \rfloor }}
            \binom{n- \sum_{i=1}^{q-1} n_i}{n_{q+1}} V_{g_q, \, n_q +a_q} V_{g_{q+1}, \, n_{q+1} + a_{q+1}}\\
    =&O\!\left(V_{g- \sum_{i= 1}^{q-1}g_i, \, n - \sum_{i= 1}^{q-1}n_i}\left(\frac{1+n^{s_q - a_q}}{g^{s_q -a_q - a_{q+1}}} + \frac{1+n^{2\lfloor\frac{g}{q+1} \rfloor - a_{q+1} }}{g^{2\lfloor\frac{g}{q+1} \rfloor -a_{q} - a_{q+1}}}\right)\right).
\end{align*}
By taking logs and comparing, for $g$ large enough for we have
\[\frac{1+n^{s_q - a_q}}{g^{s_q -a_q - a_{q+1}}} \geq \frac{1+n^{2\lfloor\frac{g}{q+1} \rfloor - a_{q+1} }}{g^{2\lfloor\frac{g}{q+1} \rfloor -a_{q} - a_{q+1}}}.\]

Thus we can bound \eqref{eq: th sommm des produits rec. séparation des q-1 premiers termes et des deux autres} by
\begin{align*}
    \dsum_{\substack{g_1 + \cdots + g_{q+1} = g \\ \bigsqcup I_i = \ens{1}{n}\\ 2g_i + \left|I_i\right| +a_i \geq s_i\\ g_{q+1} \geq \lfloor \frac{g}{q+1} \rfloor}}
            \frac{\prod_{i=1}^{q+1}V_{g_i, \, \left|I_i\right| +a_i}}{V_{g, \, n+m}}
   \leq& C \frac{1+n^{s_q - a_q}}{g^{s_q -a_q - a_{q+1}}}
            \dsum_{\substack{n_1 + \cdots + n_{q} = n \\ g_1 + \cdots + g_{q} = g\\ 2g_i + n_i +a_i \geq s_i \\ 2g_q +n_q \geq \lfloor \frac{g}{q+1}\rfloor}}
            \frac{n!}{n_1! \cdots n_q!}\frac{V_{g_q, \, n_q}\prod_{i=1}^{q-1}V_{g_i, \, n_i +a_i}}{V_{g, \, n+m}}
\end{align*}
Now, we apply the induction hypothesis to the right-hand summation. Note that on the right-hand side we have $s_{q} = \lfloor g/(q+1)\rfloor$ and $a_q =0$ and thus it is bounded by 
\[\sum_{i=1}^{q-1}\frac{1+n^{\sum_{\substack{j=1\\j \neq i}}^{q-1} s_j - a_j}}{g^{m - a_i + \sum_{\substack{j=1\\j \neq i}}^{q-1}s_j - a_j}} \leq q\frac{ 1+ n^{\sum_{\substack{j=1}}^{q-1} s_j - a_j}}{g^{m + \sum_{\substack{j=1}}^{q-1}s_j - a_j}},\]
for $g$ sufficiently large (which can be seen by taking logs and comparing as before).
So \eqref{eq: th sommm des produits rec. séparation des q-1 premiers termes et des deux autres} is bounded by
\begin{align*}
    \dsum_{\substack{g_1 + \cdots + g_{q+1} = g \\ \bigsqcup I_i = \ens{1}{n}\\ 2g_i + \left|I_i\right| +a_i \geq s_i\\ g_{q+1} \geq \lfloor \frac{g}{q+1} \rfloor}}
            \frac{\prod_{i=1}^{q+1}V_{g_i, \, \left|I_i\right| +a_i}}{V_{g, \, n+m}}
   \leq& C \frac{1+n^{s_q - a_q}}{g^{s_q -a_q - a_{q+1}}}
            \frac{ 1+ n^{\sum_{\substack{j=1}}^{q-1} s_j - a_j}}{g^{m + \sum_{\substack{j=1}}^{q-1}s_j - a_j}}
    \leq C \frac{ 1+ n^{\sum_{\substack{j=1}}^{q} s_j - a_j}}{g^{m - a_{q+1} + \sum_{\substack{j=1}}^{q}s_j - a_j}},
\end{align*}
where the $C$ changes between inequalities.
Combining this with \eqref{eq: th somme des produits rec. Decomposition de la somme principal} and the symmetry over the $j$ gives the result. 
\end{proof}

\begin{remark}
    \label{remark: sur th somme des produits}
    We make some comments on Theorem \ref{th: somme d'un produit quelconque de V_gini}:
    \begin{itemize}
        \item In the statement of the theorem we only specified that $s_i \geq 3$ because $V_{g_i, n_i +a_i} = 0$ if $2g_i + n_i +a_i < 3$. For the same reason, we can always replace $s_i$ by $\max\{s_i,a_i\}$ and so we can assume that $s_i\geq a_i$.
        \item The bound that we obtain is sharp in the sense that we can find terms in the summation whose leading order expansion matches each of the terms in the summation bound. 

        Indeed, for simplicity let us assume that $n\to\infty$ as $g\to\infty$. Then, to obtain the $i$th summand in the bound, we choose $g_j=0$ and $n_j=s_j-a_j$ for $j\neq i$. Note that we are assuming as per the previous remark that $s_j\geq a_j$ and that $n$ is sufficiently large such that $n\gg \sum_{j\neq i}s_j-a_j$. We have $g_i=g$ and $n_i=n-\sum_{j\neq i}s_j-a_j$ and so there are
            $$\frac{n!}{(n-\sum_{j\neq i} s_j-a_j)!\prod_{j\neq i}(s_j-a_j)!}=n^{\sum_{j\neq i}s_j-a_j}\left(1+O\left(\frac{1}{n}\right)\right)$$
        terms in the volume sum of this form, corresponding to the selection of cusp labels. 
        
        Since the $s_j$ are fixed, the volumes $V_{0,s_j}$ are of constant order and so the $V_{g_i,n_i+a_i}$ terms contribute an order of
            $$n^{\sum_{j\neq i}s_j-a_j}\frac{V_{g,n+a_i-\sum_{j\neq i}s_j-a_j}}{V_{g,n+m}}\asymp \frac{n^{\sum_{j\neq i}s_j-a_j}}{g^{m-a_i+\sum_{j\neq i}s_j-a_j}},$$
        by using Theorem \ref{th: developpement premier ordre comparaison classique}.
    \end{itemize}
\end{remark}

\subsection{First order expansions for moduli space volumes}
This subsection is dedicated to the proof of Theorem \ref{th: developpement premier ordre comparaison classique}. The structure of the proof relies on using the recursion relations for intersection numbers, namely Lemma \ref{lemma: recursion_rels}, as well as the general inductive proof schematic by Mirzakhani and Zograf \cite{MirzakhaniZograf2024}.

To this end, we will use bounds on the quantities appearing in Theorem \ref{th: developpement premier ordre comparaison classique}, this is Lemma \ref{lemma 5.1 MZ} below. These bounds are actually stated in Lemma 5.1 in the article of Mirzakhani and Zograf \cite{MirzakhaniZograf2024} however, some of the arguments in their proof require some more care to obtain the desired level of uniformity in $n$ and $g$. This was remedied under the assumption $\frac{n^2}{g}\to 0$ by Lipnowski and Wright \cite[Appendix A]{LipnowskiWright2024} and very similar care is needed in \cite[Section 4]{HideMaceraThomas2025}. We therefore state the bounds and divert their proof to Appendix \ref{Appendix: proof lemma 5.1}.

\begin{lemma}[{\cite[Lemma 5.1]{MirzakhaniZograf2024}}]
    \label{lemma 5.1 MZ}
   Assume that \[\lim_{g\to \infty} \frac{n}{g} = 0.\]
   Then there exist universal constants $c_0, c_1, c_2, c_3 > 0$ such that:
\begin{enumerate}[(i)]
    \item For all $k \geq 1$,
    \[
        c_0 \cdot \frac{n\mathds{1}_{2 \leq n\leq g}}{2g - 2 + n} \;\leq\; 1 - \frac{[\tau_k \tau_0^{n-1}]_{g,n}}{V_{g,n}} \;\leq\; c_1 \cdot \frac{n k^2}{2g - 2 + n}.
    \]
    \item
    \[
        \left| \frac{\left( 2g - 2 + n \right)V_{g,n}}{V_{g,n+1}} - \frac{1}{4\pi^2} \right|\;\leq\; c_2 \cdot \frac{n}{2g - 2 + n}.
    \]
    \item
    \[
        \frac{V_{g-1,n+4}}{V_{g,n+2}} \;\leq\; 1 - c_3 \cdot \frac{n}{2g - 2 + n}.
    \]
\end{enumerate}
\end{lemma}
\begin{remark}
\label{remark: on lemma 5.1 MZ}
    Part $(i)$ implies under the same assumptions on $n$ and $g$ that
    \[\left|1 - \frac{[\tau_k \tau_0^{n-1}]_{g,n}}{V_{g,n}}\right| = O\!\left(\frac{k^2(1+n)}{g}\right).\]
    In fact, in Appendix \ref{Appendix: proof lemma 5.1}, we will prove the following stronger statement: let $n_0 \in \N$ and $d= \tuple{d_1}{d_{n_0}} \in \Z_{\geq 0}^{n_0}$. Then,
    \[
        c_0 \cdot \frac{n\mathds{1}_{2 \leq n\leq g}}{2g - 2 + n} \;\leq\; 1 - \frac{[\prod_{i=1}^{n_0}\tau_{d_i} \tau_0^{n-n_0}]_{g,n}}{V_{g,n}} \;\leq\; c_1 \cdot \frac{n |d|^2}{2g - 2 + n},
    \] 
    and 
    \[\left|1 - \frac{[\prod_{i=1}^{n_0}\tau_{d_i} \tau_0^{n-n_0}]_{g,n}}{V_{g,n}}\right| = O\!\left(\frac{|d|^2(1+n)}{g}\right).\]
\end{remark}
We first use Lemma \ref{lemma 5.1 MZ} and Remark \ref{remark: on lemma 5.1 MZ} to obtain a first order expansion of the difference
\[\frac{\left[\tau_{d_1}\cdots\tau_{d_{n_0}}\tau_{0}^{n-n_0}\right]_{g,n} - \left[\tau_{d_1+1}\cdots\tau_{d_{n_0}}\tau_{0}^{n-n_0}\right]_{g,n}}{V_{g,n}}.\]

\begin{lemma}
    \label{lemma expansion 1 er ordre dif intersection number}
     Suppose that
     \[\lim_{g\to \infty} \frac{n}{g} = 0.\]
     Let $n_0 \in \Z_{\geq 0}$ and $\tuple{d_1}{d_{n_0}} \in \Z_{\geq 0}^{n_0}$, such that $n_0\leq n(g)$ for all $g$ sufficiently large.
     Then, 
     \begin{equation}
        \begin{aligned}
        \label{eq: prop du lemma expansion 1 er ordre dif intersection number}
     &\frac{\left[\tau_{d_1}\cdots\tau_{d_{n_0}}\tau_{0}^{n-n_0}\right]_{g,n} - \left[\tau_{d_1+1}\cdots\tau_{d_{n_0}}\tau_{0}^{n-n_0}\right]_{g,n}}{V_{g,n}} \\
     &= \frac{1}{2\pi^2g}\left(2(|d|-d_1) +(n-1) + \sum_{j=2}^n(2d_j +1)\delta_{1,j}\right)+  O\!\left(|d|^2\left(\frac{1}{g}+\frac{n^2}{g^2}\right)\right),
        \end{aligned} 
    \end{equation}
where, for $j$ between $2$ and $n$, $\delta_{1,j}$ is defined by
\[\delta_{1,j} = 
\begin{cases}
    0 & \text{if}\ d_1 = d_j =0,\\
    1 & \text{otherwise.}
\end{cases}
\]
\end{lemma}
\begin{proof}
    The proof is based on the third recursive formula \eqref{eq: rec III} in Lemma \ref{lemma: recursion_rels}. We have 
    \begin{align*}
   \left[\tau_{d_1}\cdots\tau_{d_{n_0}}\tau_{0}^{n-n_0}\right]_{g,n} - \left[\tau_{d_1+1}\cdots\tau_{d_{n_0}}\tau_{0}^{n-n_0}\right]_{g,n}
   =& \tilde{A}_{d,g,n} + \tilde{B}_{d, g,n} + \tilde{C}_{d,g,n},
   \end{align*}
where $\tilde{A}_{d,g,n},\, \tilde{B}_{d, g,n}\ \text{and} \ \tilde{C}_{d,g,n}$ are the differences of the $A, \, B \ \text{and},\ C$ terms in the formula \eqref{eq: rec III}, namely,
    \begin{align*}
        \tilde{A}_{d,g,n} &= 8\sum_{j=2}^n\sum_{l=0}^{d_0} (2d_j+1)(a_l-a_{l-1})\left[\tau_{d_1+d_j+l-1}\prod_{i\neq 1,j}\tau_{d_i}\right]_{g,n-1},\\
       \tilde{B}_{d, g,n}&=16\sum_{l=0}^{d_0}\sum_{k_1+k_2=l+d_1-2}(a_l-a_{l-1})\left[\tau_{k_1}\tau_{k_2}\prod_{i=2}^n\tau_{d_i}\right]_{g-1,n+1}, \\
       \tilde{C}_{d,g,n}&=16\sum_{\substack{g_1+g_2=g \\ I\sqcup J=\{2,\ldots,n\}}}\sum_{l=0}^{d_0}\sum_{k_1+k_2=l+d_1-2} (a_l-a_{l-1})\left[\tau_{k_1}\prod_{i\in I}\tau_{d_i}\right]_{g_1,|I|+1}\left[\tau_{k_2}\prod_{j\in J}\tau_{d_j}\right]_{g_2,|J|+1},       
    \end{align*}
where throughout, if $k>n_0$ then $d_k=0$.
Note that since $|d|$ is fixed, $d_0 = 3g+n-3 -|d|\to\infty$ as $g\to\infty$ and so, by Lemma \ref{lemme: aproxiamtion des coefficient a_i} (ii), 
\begin{equation}
    \label{eq: reste de la série des a_i+1 - a_i}
\sum_{i=d_0+1}^{\infty} (a_{i+1} -a_i) = O\! \left(\frac{1}{2^{6g+2n}}\right).
\end{equation}
We will see that the leading order term is given by $\tilde{A}_{d,g,n}$. We start by analysing $\tilde{B}_{d, g,n}$ by first writing
\[\frac{\tilde{B}_{d, g,n}}{V_{g,n}} = \frac{\tilde{B}_{d, g,n}}{V_{g-1,n+1}}\frac{V_{g-1,n+1}}{V_{g-1,n+2}}\frac{V_{g-1,n+2}}{V_{g,n}}.\]
Then, using Theorem \ref{Th: theorem 1.1 M_2007_1} and Lemma \ref{lemme:  aproxiamtion des coefficient a_i},
\begin{align*}
\frac{\tilde{B}_{d, g,n}}{V_{g-1,n+1}} \leq C\sum_{l=0}^{\infty}(a_{l}-a_{l-1}) \sum_{k_1 +k_2 =d_1 + l-2} 1\leq Cd_1,
\end{align*}
for some absolute $C>0$. Moreover, by Lemma \ref{Lemma 3.2 Mirzakhani2013},
$$\frac{V_{g-1,n+1}}{V_{g-1,n+2}}\frac{V_{g-1,n+2}}{V_{g,n}}\leq \frac{Cd_1}{g+n},$$
and so
\begin{equation}
    \label{eq: expansion 1 B}
\frac{\tilde{B}_{d, g,n}}{V_{g,n}} \leq \frac{Cd_1}{g}.
\end{equation}
Similarly, for $\tilde{C}_{d,g,n}$, we can use Lemma \ref{lemme initialisation de rec produit des volumes} with $a_1=a_2=1$, $s_1=s_2=3$ and $m=1$ and Lemma \ref{lemme:  aproxiamtion des coefficient a_i} to bound
\begin{align}
    \label{eq: expansion 1 C}
\frac{\tilde{C}_{d,g,n}}{V_{g,n}} &\leq 16\sum_{l=0}^{d_0}\sum_{k_1+k_2=l+d_1-2}\sum_{\substack{g_1+g_2=g \\ I\sqcup J=\{2,\ldots,n\}}}\frac{V_{g_1,|I|+1}V_{g_2,|J|+1}}{V_{g,n}}
\leq Cd_1\frac{1+n^2}{g^2}.
\end{align}
We now consider the contribution of $\tilde{A}_{d,g,n}$.
We have
\begin{equation}
    \label{eq: terme A/Vg,n-1}
    \begin{aligned}
    \frac{\tilde{A}_{d,g,n}}{V_{g,n-1}} 
    &= 8 \sum_{j=2}^{n}(2d_j +1)\left(\sum_{l=0}^{d_0-1} (a_{l+1}-a_l)\frac{\left[\tau_{d_1+d_j+l}\prod_{i\neq 1,\, j}\tau_{d_i}\right]_{g,n-1}}{V_{g,n-1}}\right. \\
    &\left.\phantom{8 \sum_{j=2}^{n}(2d_j +1)\left(\right)}+ a_0\frac{\left[\tau_{d_1+d_j-1}\prod_{i\neq 1,\, j}\tau_{d_i}\right]_{g,n-1}}{V_{g,n-1}}\right),
    \end{aligned}
\end{equation}
where, for $k > n_0$, $d_{k} =0$.
Note that
\[\left[\tau_{d_1+d_j-1}\prod_{i\neq 1,\, j}\tau_{d_i}\right]_{g,n-1}=0\]
if $d_1=d_j = 0$, and so using Remark \ref{remark: on lemma 5.1 MZ} and the tail estimate \eqref{eq: reste de la série des a_i+1 - a_i} we see that \eqref{eq: terme A/Vg,n-1} is given by
\begin{align*}
   \frac{\tilde{A}_{d,g,n}}{V_{g,n-1}} = 8 \sum_{j=2}^{n}(2d_j +1)\left(\sum_{l=0}^{+\infty} (a_{l+1}-a_l) + a_0 \delta_{1,j}\right)\left(1+ O\!\left(\frac{|d|^2(1+n)}{g}\right)\right)   + O\!\left(\frac{1}{2^{6g+2n}}\right).
\end{align*}
By Lemma \ref{lemma 5.1 MZ} (ii) we have
    $$\frac{V_{g,n-1}}{V_{g,n}}=\frac{1}{4\pi^2(2g+n-3)}\left(1+O\left(\frac{|d|^2(1+n)}{g}\right)\right).$$
Thus, since $a_0=\frac{1}{2}$ and the series of the differences of the $a_l$ equals $\frac{1}{2}$ by Lemma \ref{lemme:  aproxiamtion des coefficient a_i}, we obtain
\begin{align}
    \label{eq: expansion 1 A}
    \frac{\tilde{A}_{d,g,n}}{V_{g,n}} 
    &= \frac{1}{\pi^2(2g+n-3)}\left(2(|d|-d_1) +(n-1) + \sum_{j=2}^n(2d_j +1)\delta_{1,j}\right)\left(1+ O\!\left(\frac{|d|^2(1+n)}{g}\right)\right) \notag \\
    &= \frac{1}{2\pi^2g}\left(2(|d|-d_1) +(n-1) + \sum_{j=2}^n(2d_j +1)\delta_{1,j}\right)\left(1+ O\!\left(\frac{|d|^2(1+n)}{g}\right)\right).
\end{align}
To turn the error into an additive error rather than a multiplicative one as in the lemma statement, we note that the leading term in the product is bounded by $C\frac{n}{g}$. Thus, combining equations \eqref{eq: expansion 1 B}, \eqref{eq: expansion 1 C} and \eqref{eq: expansion 1 A} gives the result.
\end{proof}
We now prove Theorem \ref{th: developpement premier ordre comparaison classique}.
\begin{proof}[Proof of Theorem \ref{th: developpement premier ordre comparaison classique}]
We start with (i). We will assume that all of the $d_i$ are non-zero otherwise we can just work with a smaller $n_0$. We can write
\begin{align*}
1 - \frac{\left[\tau_{d_1}\cdots\tau_{d_{n_0}}\tau_{0}^{n-n_0}\right]_{g,n}}{V_{g,n}}    =& \underbrace{\sum_{k=0}^{d_1-1}\frac{\left[\tau_k\tau_0^{n-1}\right]_{g,n}-\left[\tau_{k+1}\tau_0^{n-1}\right]_{g,n}}{V_{g,n}}}_{(1)}\\
      &+\underbrace{\sum_{i=1}^{n_0-1}\sum_{k=0}^{d_{i+1}-1}\frac{\left[\tau_{k}\tau_{d_1}\cdots\tau_{d_i}\tau_{0}^{n-(i+1)}\right]_{g,n}-\left[\tau_{k+1}\tau_{d_1}\cdots\tau_{d_i}\tau_{0}^{n-(i+1) }\right]_{g,n}}{V_{g,n}}}_{(2)}.
\end{align*}
By applying Lemma \ref{lemma expansion 1 er ordre dif intersection number} the $k=0$ term in (1) evaluates to 
    $$\frac{n}{2\pi^2g}+O\left(|d|^2\left(\frac{1}{g}+\frac{n^2}{g^2}\right)\right).$$
Similarly the remaining terms in (1) evaluate to
    $$\frac{n(d_1-1)}{\pi^2g}+O\left(|d|^2\left(\frac{1}{g}+\frac{n^2}{g^2}\right)\right),$$
and so combining these terms we have
    $$(1)=\frac{(2d_1-1)n}{2\pi^2g}+O\left(|d|^2\left(\frac{1}{g}+\frac{n^2}{g^2}\right)\right).$$
For (2), we also separate out the $k=0$ case and use Lemma \ref{lemma expansion 1 er ordre dif intersection number} to see that it evaluates to
    \begin{align*}
        \frac{1}{2\pi^2g}&\left((n_0-1)(n-1)+\sum_{i=1}^{n_0-1}4(d_1+\ldots+d_i)+i\right)+O\left(|d|^4\left(\frac{1}{g}+\frac{n^2}{g^2}\right)\right)\\
        &=\frac{n(n_0-1)}{2\pi^2g}+O\left(|d|^4\left(\frac{1}{g}+\frac{n^2}{g^2}\right)\right),
    \end{align*}
where we use the fact that the part of the summation involving only the $d_i$ and $i$ is uniformly bounded by $|d|^2$ and hence lives inside the error term.
The $k\geq 1$ terms are dealt with similarly using Lemma \ref{lemma expansion 1 er ordre dif intersection number} to obtain
\begin{align*}
    \frac{1}{2\pi^2g}&\sum_{i=1}^{n_0-1}\sum_{k=1}^{d_{i+1}-1}\left(4(d_1+\ldots+d_i)+2(n-1)\right)+O\left(|d|^4\left(\frac{1}{g}+\frac{n^2}{g^2}\right)\right)\\
    &=\frac{n\sum_{i=1}^{n_0-1}(d_{i+1}-1)}{\pi^2g}+O\left(|d|^4\left(\frac{1}{g}+\frac{n^2}{g^2}\right)\right),
\end{align*}
where again we have used that the part of the summation involving only the $d_i$ can be absorbed into the error term.

Combining these together results in 
\begin{align*}
    (2) = \frac{n\sum_{i=2}^{n_0}d_i}{\pi^2g}-\frac{n(n_0-1)}{2\pi^2g}+O\left(|d|^4\left(\frac{1}{g}+\frac{n^2}{g^2}\right)\right).
\end{align*}
So summing (1) and (2) we obtain 
\begin{equation}
    \label{eq: expansion 2ème ordre intersection number sur volume.}
    1 - \frac{\left[\tau_{d_1}\cdots\tau_{d_{n_0}}\tau_{0}^{n-n_0}\right]_{g,n}}{V_{g,n}}
    =\frac{\left(2|d|-n_0\right)n}{2\pi^2g} + O\!\left(|d|^4\left(\frac{1}{g}+\frac{n^2}{g^2}\right)\right),
\end{equation}
as desired.

We now prove (ii) using recursive formula Lemma \ref{lemma: recursion_rels} \eqref{eq: rec II}. We have
\begin{equation*}
\frac{(2g+n-2)V_{g,n}}{V_{g,n+1}} 
    = \frac{1}{2}\sum_{l=0}^{3g+n-3}(-1)^l(l+1)\frac{\pi^{2l}}{(2l+3)!}\frac{\left[\tau_{l+1}\tau_{0}^{n}\right]_{g,n+1}}{V_{g,n+1}}.
\end{equation*}
Since
    $$\sum_{l=0}^{+ \infty} (-1)^l(l+1)\frac{\pi^{2l}}{(2l+3)!} = \frac{1}{2\pi^2},$$
we obtain
    \begin{align}
    \label{eq: first step (ii) th expansion form}
        \frac{(2g+n-2)V_{g,n}}{V_{g,n+1}} - \frac{1}{4\pi^2}
    = \frac{1}{2}\sum_{l=0}^{3g+n-3}(-1)^l(l+1)\frac{\pi^{2l}}{(2l+3)!}\left(\frac{\left[\tau_{l+1}\tau_{0}^{n}\right]_{g,n+1}}{V_{g,n+1}}-1\right)+O\left(\frac{1}{2^{3g}}\right),
    \end{align}
where the error term arises from bounding (using Stirling approximation) 
    $$\frac{(l+1)\pi^{2l}}{(2l+3)!}\leq C\frac{(l+1)\pi^{2l}}{(2l+3)^{2l+3}}\leq C\frac{1}{l^2}\frac{1}{2^{3g+n}},$$
for $l\geq 3g+n-2$ and $g$ sufficiently large resulting in
    $$\sum_{l=3g+n-2}^{\infty} (-1)^l(l+1)\frac{\pi^{2l}}{(2l+3)!}=O\left(\frac{1}{2^{3g}}\right).$$
Now, continuing on \eqref{eq: first step (ii) th expansion form} and using \eqref{eq: expansion 2ème ordre intersection number sur volume.}, we obtain 
    \begin{align*}
        \frac{(2g+n-2)V_{g,n}}{V_{g,n+1}} - \frac{1}{4\pi^2}
    &= -\frac{1}{2}\sum_{l=0}^{3g+n-3}(-1)^l(l+1)\frac{\pi^{2l}}{(2l+3)!}\left(\frac{(2(l+1)-1)(n+1)}{2\pi^2g}+O\left(\ell^2\left(\frac{1}{g}+\frac{n^2}{g^2}\right)\right)\right)\\
    &=-\frac{n+1}{2\pi^2g}\sum_{l=0}^{3g+n-3}(-1)^l(l+1)^2\frac{\pi^{2l}}{(2l+3)!}\\
    &\ \ +\frac{n+1}{4\pi^2g}\sum_{l=0}^{3g+n-3}(-1)^l(l+1)\frac{\pi^{2l}}{(2l+3)!}+O\left(\ell^2\left(\frac{1}{g}+\frac{n^2}{g^2}\right)\right).
    \end{align*}
Since 
    $$\sum_{l=0}^{ \infty} (-1)^l(l+1)^2\frac{\pi^{2l}}{(2l+3)!} = -\frac{1}{4\pi^2},$$
and the tails of these series are $O(2^{-3g})$ as before, we obtain the desired
\[\frac{(2g+n-2)V_{g,n}}{V_{g,n+1}} - \frac{1}{4\pi^2} = \frac{n}{4\pi^4g} + O\!\left(\frac{1}{g}+\frac{n^2}{g^2}\right).\]
We now turn to proving (iii) using the formula Lemma \ref{lemma: recursion_rels} \eqref{eq: rec I} to write
\begin{align*}
    \frac{V_{g-1,n+4}}{V_{g,n+2}}=\frac{\left[\tau_{1}\tau_{0}^{n+1}\right]_{g,n+2}}{V_{g,n+2}} - 6\sum_{\substack{g_1+g_2 =g \\ I \bigsqcup J = \ens{1}{n}}} \frac{V_{g_1, |I|+2}V_{g_2, |J|+2}}{V_{g,n+2}}.
\end{align*}
We distinguish the two terms in the summation when $2g_1+|I|+2=3$ and $2g_2+|J|+2$ to write

\begin{align*}
    \frac{V_{g-1,n+4}}{V_{g,n+2}}=\frac{\left[\tau_{1}\tau_{0}^{n+1}\right]_{g,n+2}}{V_{g,n+2}}- 12n \frac{V_{g,n+1}}{V_{g,n+2}} - 6\sum_{\substack{g_1+g_2 =g \\ I \bigsqcup J = \ens{1}{n} \\ 2g_1 + |I| + 2 \geq 4 \\ 2g_2 + |J| + 2 \geq 4}} \frac{V_{g_1, |I|+2}V_{g_2, |J|+2}}{V_{g,n+2}}.
\end{align*}
We now use \eqref{eq: expansion 2ème ordre intersection number sur volume.} for the first term, Lemma \ref{lemma 5.1 MZ} (ii) for the second and Lemma \ref{lemme initialisation de rec produit des volumes}  with $a_1=a_2=2$, $s_1=s_2=4$ and $m=2$ for the third term to give
\[\frac{V_{g-1,n+4}}{V_{g,n+2}} =  1- \frac{2n}{\pi^2g} + O\!\left(\frac{1}{g}+\frac{n^2}{g^2}\right).\]
\end{proof}
We conclude this section with the corresponding expansion for bordered volumes.
\begin{lemma}
   \label{lemma: dl de Vgn(x)/Vgn}
   Suppose that
   \[ \lim_{g \to \infty} \frac{n}{g} = 0.\]
   Let $K \in \N$ and $(x_1, \, \cdots, \, x_K) \in \R_{>0}^K$. 
   Then,
   \[\frac{V_{g,n}\tuple{2x_1}{2x_K}}{V_{g,n}} = \prod_{i=1}^{K}\frac{\sinh(x_i)}{x_i}\left(1 + O\!\left(\frac{\left(\sum_{i=1}^{K}x_i^2\right)(1+n)}{g}\right)\right)\]
\end{lemma}
\begin{proof}
Using Theorem \ref{Th: theorem 1.1 M_2007_1} we write
\begin{align*}
   \frac{V_{g,n}\tuple{2x_1}{2x_K}}{V_{g,n}}
   =& \sum_{|d| \leq 3g - 3 + n} \left(\frac{\left[ \tau_{d_1}, \dots, \tau_{d_K}, \tau_0^{n-K} \right]_{g,n}}{V_{g,n}}-1+1\right)\prod_{i=1}^K \frac{x_i^{2d_i}}{(2d_i + 1)!}.\\
\end{align*}
By Remark \ref{remark: on lemma 5.1 MZ}, 
    \begin{align*}
        \frac{\left[ \tau_{d_1}, \dots, \tau_{d_K}, \tau_0^{n-K} \right]_{g,n}}{V_{g,n}}-1=O\left(|d|^2\left(\frac{1+n}{g}\right)\right)
    \end{align*}
so that 
\begin{align*}
   \frac{V_{g,n}\tuple{2x_1}{2x_K}}{V_{g,n}}
   &= \sum_{|d| \leq 3g - 3 + n} \left(1+O\left(|d|^2\frac{1+n}{g}\right)\right)\prod_{i=1}^K \frac{x_i^{2d_i}}{(2d_i + 1)!}\\
   &=\prod_{i=1}^K\frac{\sinh(x_i)}{x_i}-\sum_{|d|>3g+n-3}\prod_{i=1}^K\frac{x_i^{2d_i}}{(2d_i+1)!}\\
   &\ \ \ \ \ +O\left(\frac{1+n}{g}\sum_{d\in\Z_{\geq 0}^K}\prod_{i=1}^K\frac{x_i^{2d_i}}{(2d_i+1)!|d|^2}\right).
\end{align*}
It is clear that the tail series will lie in the error term and so we focus on the big $O$ term. We have
\begin{align*}
   \sum_{d =\tuple{d_1}{d_K} \in \N^K}|d|^2\prod_{i=1}^K \frac{x_i^{2d_i}}{(2d_i + 1)!}
   &\leq C\sum_{d =\tuple{d_1}{d_K} \in \N^K} \sum_{i=1}^{K}d_i^2 \prod_{j=1}^K\frac{x_j^{2d_j}}{(2d_j + 1)!}\\
   &\leq C\sum_{i=1}^{K}\prod_{j \neq i} \frac{\sinh(x_j)}{x_j}\sum_{d_i = 1}^{\infty}d_i^2\frac{x_i^{2d_i}}{(2d_i + 1)!}\\
   &\leq C\sum_{i=1}^{K}\prod_{j \neq i} \frac{\sinh(x_j)}{x_j}\sum_{d_i = 1}^{\infty}\frac{x_i^{2d_i}}{(2d_i -1)!}\\
   &\leq C\sum_{i=1}^{K}x_i^2 \prod_{j=1}^{K} \frac{\sinh(x_j)}{x_j},
\end{align*}
and so put together, we obtain the expansion.
\end{proof}
\section{Poisson statistics}
\label{sec: Poisson statistic}
This section is dedicated to proving Theorems \ref{th: stat de poisson dans le cas n^2/g -> 0}, \ref{th: stat de poisson dans le cas n^2/h -> +infty} and \ref{thm:stat de poisson dans le cas n^2/g -> alpha}. The main strategy of proof in each case will be to understand the dominant contribution to the factorial moments. To this end, we make the following definition. 
\begin{definition}
    Let $k\in\N$ and $q\in\Z_{\geq 0}$. A multicurve $\Gamma = \tuple{\gamma_1}{\gamma_k}$ on $\Sigma_{g,n}$ is called a $q$-pants multicurve of length k if 
        \begin{itemize}
            \item exactly $q$ of the curves $\gamma_{i_1},\ldots,\gamma_{i_q}$ are separating and cut the surface into 
                $$\Sigma_{g,n}\setminus\gamma_{i_p}=\Sigma_{0,3}\sqcup\Sigma_{g,n-1},$$
            \item the remaining $k-q$ curves are non-separating and together do not separate the surface. In other words, 
                $$\Sigma_{g,n}\setminus\Gamma = \Sigma_{g-(k-q),n-2k-3q}\sqcup\bigsqcup_{i=1}^q\Sigma_{0,3}$$
        \end{itemize}
    We call $\Gamma$ a full pants multicurve of length $k$ if $q=k$.
\end{definition}

\begin{remark}
    A full pants multicurve is an example of what is termed an unnested multicurve in \cite{HideThomas2024}.
\end{remark}

We begin by counting the number of mapping class group orbits of $q$-pants multicurves of length $K$.
\begin{lemma}
    \label{lemma: nombre de pants multicurve}
    Let $K\in\N$ and $q\in\Z_{\geq 0}$ such that $q\leq K$. The number of distinct mapping class group orbits of q-pants multicurves of length $K$ on $\Sigma_{g,n}$ is
    \[\binom{K}{q}\frac{n!}{(n-2q)!}.\]
    Furthermore, the positions of the $q$ separating curves is invariant for an orbit and upon fixing these positions, the number of mapping class group orbits is 
    $$\frac{n!}{(n-2q)!}.$$
\end{lemma}
To prove Lemma \ref{lemma: nombre de pants multicurve} we begin with the full pants multicurves of length $1$.
\begin{lemma}
    \label{lemma: nombre d'orbite de pants}
    Let $g,n\in\Z_{\geq 0}$ satisfy $2g+n \geq 3$ and $n\geq 2$. Then the number of mapping class group orbits of full pants multicurves of length 1 in $\Sigma_{g,n}$ is 
    \[\frac{n!}{(n-2)!}.\]
\end{lemma}
\begin{proof}
    Let $\gamma$ be a full pants multicurve of length 1. Since the mapping class group fixes cusps, if $\gamma$ separates off the cusps $\{i,j\}$ and $\tilde{\gamma}$ separates off the cusps $\{\tilde{i},\tilde{j}\}$ then $\gamma$ and $\tilde{\gamma}$ can only be in the same mapping class group orbit if $\{i,j\}=\{\tilde{i},\tilde{j}\}$. Thus we have at minimum  
    \begin{equation*}
    \binom{n}{2}    
    \end{equation*}
    different orbits.
    
    We now prove that for every pair $\{i,j\}$ of cusps, we have exactly two different orbits. First remark that there is a unique primitive, simple closed unoriented geodesic that separates the pants containing two given cusps and so we just have to show that when we orient the geodesic, that it is in a different mapping class group orbit to its reverse orientation.
    We proceed by contradiction. Assume that $\gamma$ and $\gamma^{-1}$ are in the same orbit so that there exists $f$ a homeomorphism of $\Sigma_{g,n}$ such that $f(\gamma) = \gamma^{-1}$. By the collar lemma \cite[Theorem 4.4.6]{buser2010geometry}), there exists an annular neighbourhood $N$ around the curve $\gamma$ and because $f$ is a homeomorphism, $f(N)$ is homeomorphic to an annulus. 
    
    Since $f(N)$ contains $\gamma$, there exists a homeomorphism $\varphi$ isotopic to the identity map on $\Sigma_{g,n}$ such that $\varphi \circ f(N) = N$. Thus we may assume that $f(N) = N$ and so in this case, 
    \[f|_{\Sigma_{g,n}\setminus N} : \Sigma_{g,n}\setminus N \xrightarrow{} \Sigma_{g,n}\setminus N \]
    is a homeomorphism. But, $\Sigma_{g,n}\setminus N$ has two connected components; one is the pair of pants with the two cusps $i$ and $j$, $P$, and the second is the rest of the surface that we denote $B$. Since $f$ fixes the cusps we must have $f(P) = P$ and $f(B) =B$ and so $f$ fixes the two boundary components of $N$ set-wise. In other words, $f|_N \in \text{MCG}(N) \simeq \Z$. However, the mapping class group of an annulus is generated by Dehn twists which do not change the orientation of $\gamma$. This gives the contradiction and completes the proof.
\end{proof}
We now prove Lemma \ref{lemma: nombre de pants multicurve}.
\begin{proof}[Proof of Lemma \ref{lemma: nombre de pants multicurve}]
     Let 
    \[\gamma_1, \cdots, \gamma_K\]
    be simple oriented closed geodesics in which exactly $q$ bound a pair of pants (with two cusps) that are pairwise distinct and the $K-q$ other curves are non-separating. We will count the number of $q$-pants multicurves in different mapping class group orbits that we can create with these curves. 
    
    Note that since the curves bounding the cusps are disjoint and distinct, they separate off different cusps and hence are in different mapping class group orbits. Moreover, since non-separating closed curves are not mapping class group equivalent to separating ones, the position of the separating curves is an invariant of the orbit. We clearly have
    \begin{equation}
        \label{eq: nb de choix pour les curves qui bound  un pant}
        \binom{K}{q}    
    \end{equation}
    different choices for the positions of the separating curves. The $K-q$ non-separating curves are in the same mapping class group orbit (see for example \cite[Section 1.3.1]{farb2011primer}) and can be interchanged one by one. Thus the number of orbits after fixing the positions of the separating curves is determined by the number of ways to choose the separating curves. \\
    \\
    The number of choices for the first curve is given by Lemma \ref{lemma: nombre de pants multicurve} as 
        $$\frac{n!}{(n-2)!}.$$
    For the second curve, now that the first two cusps are chosen, has 
        $$\frac{(n-2)!}{(n-4)!},$$
    choices and so on which, upon multiplying together, give a factor
        $$\frac{n!}{(n-2q)!}.$$
    This gives the second claim. The first claim comes from multiplying this by the separating curve position choices
    \[ \binom{K}{q}\frac{n!}{(n-2q)!}.\]

\end{proof}

  We remark that for the purposes of the Mirzakhani integration formula, we consider oriented multicurves. However, for our theorems, we wish to count the unoriented curves. As a result we must multiply through by a pre-factor of $2^{-K}$ where $K$ is the multicurve length. Now recall the definition of the random variable $Y_{g,n,r_1,\ldots,r_k}$ from Section \ref{subsection : Factorial moments and particular notations}.

\subsection{Proof of Theorem \ref{th: stat de poisson dans le cas n^2/g -> 0}}
\label{subsection : stat poisson n^2/g ->0}
In this section, we will prove Theorem \ref{th: stat de poisson dans le cas n^2/g -> 0} using Mirzakhani's integration formula and the method of factorial moments. The leading order contribution to the product of the factorial moments will arise from the multicurves that are simple and non-separating, this is analogous to the closed surface case.

Recall the decomposition 
    \[(N_{g,n, [a_1,b_1]})_{r_1} \cdots (N_{g,n, [a_k,b_k]})_{r_k}= Y_{g,n,r_1, \cdots, r_k}+Y_{g,n,r_1, \cdots, r_k}'\]
from Section \ref{subsection : Factorial moments and particular notations}.
\begin{lemma}
    \label{lemma: n^2/g -> 0, multicurve}
    Let $k,r_1,\,\cdots,\, r_k\in\N$ 
        $$0\leq a_1<b_1< a_2< \cdots <b_{k-1}< a_k <b_k.$$
    Assume that 
    \[\frac{n^2}{g} \xrightarrow[g \to \infty]{} 0.\] 
    Then,
    \[\EE\left(Y_{g,n,r_1, \cdots, r_k}\right) = \prod_{i=1}^{k}\lambda_{i}^{r_i} + O\!\left(\frac{1+n^2}{g}\right),\]
    where 
    \[\lambda_i = \int_{a_i}^{b_i}\frac{2\sinh(\frac{x}{2})^2}{x}\, \df x.\]
\end{lemma}
\begin{proof}
    Let $K = \sum_{i=1}^{k}r_i$ and $A = \prod_{i=1}^{k}[a_i, b_i]^{r_i}$.
    By Mirzakhani's integration formula Lemma \ref{lemma: Mirzakhani integration formula}, the expected value of $Y_{g,n,r_1, \cdots, r_k}$ is given by 
    \begin{equation}
        \label{eq: decompostion en fct de nb de composante connexe de E(Y_g,n,r_1,, r_k) 2}
        \EE\left(Y_{g,n,r_1, \cdots, r_k}\right) = \sum_{q=1}^{K+1} \sum_{\substack{[\Gamma], \, \text{s.t.}\\q(\Gamma) = q}} 2^{-M(\Gamma)-K}\int_{A}\frac{V_{g,n}(\Gamma,x)}{V_{g,n}}x_1\cdots x_K \df x_1\wedge \cdots \wedge \df x_K,
    \end{equation}
where $q(\Gamma)$ denotes the number of connected components of $\Sigma_{g,n}\setminus\Gamma$.
We first compute the contribution of $q=1$. In this case, by Lemma \ref{lemma: nombre de pants multicurve} there is a single mapping class group orbit. For any multicurve $\Gamma_0$ in this orbit, we have $\Sigma_{g,n}\setminus\Gamma_0= \Sigma_{g-K,n+2K}$ and $M(\Gamma_0)=0$. Thus the contribution is given by
\begin{equation}
    \label{eq: n^2/g-> 0 cas q=1}
    \begin{aligned}
  &\frac{1}{2^K} \frac{V_{g-K,n+2K}}{V_{g,n}}\int_{A}\frac{V_{g-K,n+2K}(x_1,x_1, \cdots, x_K,x_K)}{V_{g-K,n+2K}}\left(\prod_{i=1}^{K}x_i\right) \bigwedge_{i=1}^{K} \df x_i \\
  &=\left(1+O\left(\frac{1+n}{g}\right)\right)\int_{A}\frac{\prod_{i=1}^K 2\sinh^2\left(\frac{x_i}{2}\right)}{\prod_{i=1}^K x_i}\left(1+O\left(\frac{\sum_{i=1}^K  x^2_i(1+n)}{g}\right)\right)\bigwedge_{i=1}^{K} \df x_i\\
  &=\prod_{i=1}^{k}\lambda_{i}^{r_i} + O\!\left(\frac{1+n}{g}\right),
    \end{aligned}
\end{equation}
where on the second line we used Lemma \ref{lemma: dl de Vgn(x)/Vgn} to expand the bordered moduli space volume and Theorem \ref{th: developpement premier ordre comparaison classique} (iii) for the volume pre-factor.

We now turn to the remaining possibilities for $q$. Let $\Gamma$ be a multicurve such that $\Sigma_{g,n} \setminus \Gamma$ has $q$ connected components. The $n$ cusps are distributed on the $q$ subsurfaces and each subsurface has at least one bordered component coming from a component of the multicurve $\Gamma$. Thus, there exist $q$ disjoint subsets
 \[\bigsqcup_{i=1}^{q} I_i = \ens{1}{n},\]
 corresponding to the labelled cusps in each subsurface and $a_1, \cdots, \, a_q\in\N$ corresponding to the number of boundaries created by the multicurve $\Gamma$ on each subsurface so that 
 \[\sum_{i=1}^q a_i = 2K.\]
 Then,
 \[\Sigma_{g,n} \setminus \Gamma = \bigsqcup_{i=1}^q \Sigma_{g_i, |I_i| + a_i}.\]
 An Euler characteristic computation shows that the genera $g_i$ in the resulting subsurfaces satisfy
 \[\sum_{i=1}^{q} g_i = g+q-K-1.\]
 Any multicurve in the same mapping class group orbit as $\Gamma$ will give rise to the same data so that the collection $\{(g_i,a_i,I_i)\}$ is well-defined for the orbit.

On the other hand, the number of mapping class group orbits giving rise to the same decomposition is bounded in terms of $K$ only. Indeed, any permutation of the curve components in the multicurve gives rise to the same decomposition and potentially a new mapping class group orbit (so $K!$ choices). Aside from this, any other orbit corresponds to a gluing of the subsurfaces to recover $\Sigma_{g,n}$ which is a matching of the boundary components. For this there are at most $(2k-1)!!$ choices (not all matchings are permissible so this is an over-count). This means there are at most $K!(2k-1)!!$ orbits for each decomposition (see also \cite[Lemma 3.4]{HideThomas2024}). It follows that
 \begin{align}
    &\sum_{q=2}^{K+1} \sum_{\substack{[\Gamma], \, \text{s.t.}\\q(\Gamma) = q}} 2^{-M(\Gamma)-K}\int_{A}\frac{V_{g,n}(\Gamma,x)}{V_{g,n}}x_1\cdots x_K \df x_1\wedge \cdots \wedge \df x_K \notag\\
&\leq C\sum_{q=2}^{K+1} \sum_{\substack{a_1 +\cdots+a_q = 2K\\ a_i \geq 1 \\ I_1\sqcup\ldots\sqcup I_q=\{1,\ldots,n\}\\ g_1 + \cdots + g_q = g+q -K-1}}
         \frac{\prod_{i=1}^{q}V_{g_i,n_i + a_i}}{V_{g,n}}
         \int_{A} \prod_{i=1}^{K}\frac{2\sinh(\frac{x_i}{2})^2}{x_i}\df x_i \notag\\
&\leq C\sum_{q=2}^{K+1}\sum_{\substack{a_1 + \cdots +a_q = 2K\\ a_i \geq 1}}
          \sum_{\substack{ I_1\sqcup\ldots\sqcup I_q=\{1,\ldots,n\}\\ g_1 + \cdots + g_q = g+q -K-1\\ 2g_i + n_i +a_i \geq \max(3, a_i)}}
          \frac{\prod_{i=1}^{q}V_{g_i,n_i + a_i}}{V_{g,n}}
          \prod_{i=1}^{K}\lambda_i^{r_i}.
\label{eq: n^2/g-> 0 cas q entre 2 et K+1}
 \end{align}
 where on the second line we used the trivial $\sinh$ bound for bordered moduli space volumes. For a fixed $q$ and tuple $\ens{a_1}{a_q}$, the interior summation is bounded using Theorem \ref{th: somme d'un produit quelconque de V_gini} and Lemma \ref{Lemma 3.2 Mirzakhani2013} by
 \begin{align*}
 \sum_{\substack{ I_1\sqcup\ldots\sqcup I_q = \{1,\ldots,n\}\\ g_1 + \cdots + g_q = g+q -K-1\\ 2g_i + n_i +a_i \geq \max(3, a_i)}}
          \frac{\prod_{i=1}^{q}V_{g_i,n_i + a_i}}{V_{g,n}} 
    &=\frac{V_{g+q-K-1,n}}{V_{g,n}}\sum_{\substack{ I_1\sqcup\ldots\sqcup I_q = \{1,\ldots,n\}\\ g_1 + \cdots + g_q = g+q -K-1\\ 2g_i + n_i +a_i \geq \max(3, a_i)}}
          \frac{\prod_{i=1}^{q}V_{g_i,n_i + a_i}}{V_{g+q-K-1,n}} \\
    &\leq C \frac{1}{g^{2K+2-2q}}\sum_{i= 1}^{q}\frac{1+n^{\sum_{\substack{j=1\\j \neq i}}^{q} \max(3, a_j) - a_j}}{g^{- a_i + \sum_{\substack{j=1\\j \neq i}}^{q}\max(3, a_j) - a_j}}\\
    & = C \frac{1}{g^{2-2q}}
      \sum_{i= 1}^{q}
      \frac{1+n^{\sum_{\substack{j=1\\j \neq i}}^{q} \max(3, a_j) - a_j}}{g^{\sum_{\substack{j=1\\j \neq i}}^{q}\max(3, a_j)}}\\
      &\leq C\left(\frac{1+n^2}{g}\right)^{q-1},
 \end{align*}
 where in the last line we have used that because $n/g\to 0$ as $g\to\infty$, the maximum of the summation is attained when all of the $a_i$ are minimal and equal to 1.
 
 Put together, and using the fact that the number of choices of the $a_i$ is bounded in terms of $K$ only, we obtain 
 \begin{equation}
    \label{eq: n^2/g-> 0 cas q entre 2 et K+1 result}
    \sum_{q=2}^{K+1} \sum_{\substack{[\Gamma], \, \text{s.t.}\\q(\Gamma) = q}} 2^{-M(\Gamma)-K}\int_{A}\frac{V_{g,n}(\Gamma,x)}{V_{g,n}}x_1\cdots x_K \df x_1\wedge \cdots \wedge \df x_K
    \leq C \left(\frac{1+n^2}{g}\right).
 \end{equation}
 The desired result is obtained by summing equations \eqref{eq: n^2/g-> 0 cas q=1} and \eqref{eq: n^2/g-> 0 cas q entre 2 et K+1 result} 
 \end{proof}
To complete the proof of Theorem \ref{th: stat de poisson dans le cas n^2/g -> 0}, we show that the remaining curves, those counted by $Y'_{g,n,r_1,\ldots,r_k}$, are contained in an error term.
\begin{lemma}
    \label{lemma: self-insection negligeable}
    We have
    \[\EE\left(Y_{g,n,r_1,\cdots,r_k}^{\prime}\right) = O\!\left(\frac{1+n^2}{g}\right).\]
\end{lemma}
Our proof is a generalisation of the proof of \cite[Proposition 4.5]{MirzakhaniPetri2021}. One of the key inputs to this is a bound on the number of oriented closed geodesics of length at most $L$ on a finite area surface that are not iterates of oriented closed geodesics of length $<2\mathrm{arcsinh}(1)$. The corresponding bound for closed surfaces is \cite[Lemma 6.6.4]{buser2010geometry} and the proof is analogous. 
\begin{lemma}[{\cite[Lemma 6.6.4]{buser2010geometry}}]
    \label{Lemma: lemma 6.6.4 Buser}
    Let S be a hyperbolic surface with genus $g$, $n$ cusps and $b$ boundaries. For any $L>0$, there are at most
    \[C(2g+n+b-2)e^{L+6}\]
    oriented closed geodesics of length at most $L$ on $S$ which are not iterates of closed geodesics of length smaller than $2\, \mathrm{arcsinh}(1)$.
\end{lemma}
\begin{proof}
    The proof is almost identical to \cite[Lemma 6.6.4]{buser2010geometry}. The only change is that the area of the surface is now $2\pi(2g+n+b-2)$, and we use the collar lemma for non-compact surfaces  \cite[Theorem 4.4.6]{buser2010geometry}.
\end{proof}
\begin{proof}[Proof of Lemma  \ref{lemma: self-insection negligeable}]
    Let $h$ be a function defined by:
    \[\begin{array}{rclc}
        h: & \R_{>0} & \longrightarrow & \R_{\geq 0} \\
        & x & \longmapsto & \log\left(\coth\left(\frac{x}{4}\right)\right)
    \end{array}.\]
    Let $\Gamma = \ens{\gamma_1}{\gamma_K}$ be a collection of curves with either at least one curve self-intersecting or some pair having non-empty intersection.
    As in the proof of \cite[Proposition 4.5]{MirzakhaniPetri2021} (which uses \cite[Theorem 1.2 and Lemma 2.2]{AraBasmajian}) it is possible to bound the total number of intersections between the curves and themselves (all counted with multiplicity) by
    \[
\sum_{i=1}^{k} 2 r_i e^{2 b_i}
\;+\;
\binom{\sum_{i=1}^{k}r_i}{2} \left( \frac{b_k}{h(b_k)+1} \right)^2
\;=:\; C.
\]
    Now consider the regular neighbourhood of the collection of curves, denote it by $N(\Gamma)$. Since there is at least one intersection, this regular neighbourhood contains a component with Euler characteristic $\leq -1$. Moreover, since the neighbourhood retracts back onto the curve collection, the total Euler characteristic of $N(\Gamma)$ (and hence this component) is bounded below by the Euler characteristic of the graph $\mathcal{G}$ associated to the curves with vertices at each intersection point and edges the arc segments between intersections.

    For a vertex $v$ in this graph, let $\deg (v)$ denote its degree. We clearly have 
        $$|E(\mathcal{G})|=\frac{1}{2}\sum_{v\in V(\mathcal{G})}\deg(v).$$
    An intersection creates at most two incoming and outgoing arcs at a point, and so contributes at most 4 to the degree of a vertex. Let $n(\Gamma)$ the number of intersection of $\Gamma$ with multiplicity, then the total number of edges in $\mathcal{G}$ is bounded by $2n(\Gamma)$ and the Euler characteristic is bounded below by $-n(\Gamma)$. Thus since the total number of intersections (with multiplicity) is at most $C$. This results in the Euler characteristic of $\mathcal{G}$ being bounded below by $-C$.

    So, from a collection of curves $\Gamma$, we can construct a multicurve using the simple closed geodesics in the free homotopy classes of the boundary components of $N(\Gamma)$ (throwing away repeats) of size at most $2K$ whose lengths sum to at most 
        $$2\sum_{i=1}^k b_ir_i,$$
    and that when cut along, they separate off a component of signature $(g',n',b')$ satisfying 
        $$1\leq 2g'+n'+b'-2\leq C.$$

    Several collections $\Gamma$ may give rise to the same multicurve. But the maximum possible number of collections is clearly bounded by the number of ways to select $K$ closed geodesics that live inside $N(\Gamma)$ with length at most $b_k$. 
    
    Since the Euler characteristic of every component of $N(\Gamma)$ is bounded in absolute value by $C$ and there are most $K+1$ components, Lemma \ref{Lemma: lemma 6.6.4 Buser} gives at most $2(K+1)(2g'+n'+b')\exp(b_k+6)$ possible unoriented closed geodesics in $N(\Gamma)$ with length at most $b_k$ that are not iterates of closed geodesics of length at most $2\mathrm{arcsinh}(1)$. Moreover, the number of closed geodesics of length at most $2\mathrm{arcsinh}(1)$ in $N(\Gamma)$ is also only bounded in terms of $g',n'$ and $b'$. We then choose $K$ of these curves.

    In summary, there is a bound $B=B(r_1,\ldots,r_k,b_1,\ldots,b_k)$ on the number of collections $\Gamma$ giving rise to the same multicurve. Clearly then
        $${Y}_{g,n,r_1,\ldots,r_k}^{\prime}\leq \sum_{\substack{(g',n',b') \\
        1\leq 2g'+n'+b'-2 \leq C\\
        b'\geq 1}}\sum_{[\Gamma]}\sum_{\Gamma'\in [\Gamma]} \1\left(\ell_{\Gamma'}(X)\leq 2\sum_{i=1}^k b_ir_i\right),$$
where $[\Gamma]$ runs over the mapping class group orbits of all multicurves of size at most $2K$ that separate off a component of signature $(g',n',b')$ in $\Sigma_{g,n}$. 

We can now apply Mirzakhani's integration formula Lemma \ref{lemma: Mirzakhani integration formula} along with trivial $\sinh$ bounds on bordered moduli space volumes to obtain 
    \begin{align*}
        &\EE\left(Y_{g,n,r_1,\ldots,r_k}^{\prime}\right)\\
        &\leq B\sum_{\substack{(g',n',b') \\
        1\leq 2g'+n'+b'-2 \leq C\\
        b'\geq 1}}\sum_{m=1}^{2K}\sum_{q=1}^{m}\sum_{\substack{a_1+\ldots+a_q=2m-b' \\ g_1+\ldots+g_q=g-g' \\ I_1\sqcup\ldots\sqcup I_{q+1}=\{1,\ldots,n\} \\
        |I_{q+1}|=n'}} \frac{V_{g',n'+b'}\prod_{i=1}^q V_{g_i,|I_i|+a_i}}{V_{g,n}}\\
        &\hspace{2cm}\cdot\int_{\sum_{i=1}^{m}x_i\leq 2\sum_{i=1}^kr_ib_i}\prod_{i=1}^m\frac{2\sinh^2\left(\frac{x_i}{2}\right)}{x_i}\bigwedge_{i=1}^m\df x_i.
    \end{align*}
Here the second summation runs over the possible number of components in the multicurve and the third summation runs over the possible number of components obtained (other than $\Sigma_{g',n'+b'}$) after cutting the surface by the multicurve.

It is clear that the integral term is bounded uniformly only in the $r_i$ and $b_i$. Moreover, using Theorem \ref{th: somme d'un produit quelconque de V_gini} (as in the proof of Lemma \ref{lemma: n^2/g -> 0, multicurve}) the interior summation is 
    $$O\left(\frac{1+n^2}{g}\right).$$
The remaining outer summations are again bounded only in terms of the $r_i$ and $b_i$ and hence we can conclude the result.
\end{proof}
Theorem \ref{th: stat de poisson dans le cas n^2/g -> 0} is now just the combination of Lemmas \ref{lemma: n^2/g -> 0, multicurve} and \ref{lemma: self-insection negligeable}.
\begin{proof}[Proof of Theorem \ref{th: stat de poisson dans le cas n^2/g -> 0}]
    Recall that
    \[(N_{g,n, [a_1,b_1]})_{r_1} \cdots (N_{g,n, [a_k,b_k]})_{r_k} = Y_{g,n,r_1,\cdots,r_k} + Y_{g,n,r_1,\cdots,r_k}^{\prime}.\]
    Then by Lemmas \ref{lemma: n^2/g -> 0, multicurve} and \ref{lemma: self-insection negligeable},
    \[\EE\left((N_{g,n, [a_1,b_1]})_{r_1} \cdots (N_{g,n, [a_k,b_k]})_{r_k}\right) = \prod_{i=1}^{k}\lambda_i +O\!\left(\frac{1+n^2}{g}\right).\]
    Taking the limit in the latter equation and using Theorem \ref{th : Method of factorial moments} provides the desired result.
\end{proof}

\subsection{Proof of Theorem \ref{th: stat de poisson dans le cas n^2/h -> +infty}}
In this subsection, we work under the assumptions:
\[\frac{n^2}{g} \xrightarrow[g\to \infty]{} \infty \qquad \text{and} \qquad \frac{n}{g} \xrightarrow[g\to \infty]{} 0 .\]
In Section \ref{subsection : stat poisson n^2/g ->0} we saw error terms in the factorial moments of the order $n^2/g$ and without much difficulty, it is not hard to see that there are terms that contribute at this order. Thus under our new assumption of $n^2/g\to\infty$, to hope for finite factorial moments, we are required to renormalise the interval lengths that we look at. 

Thus, let $k\in\N$ and 
    $$0\leq a_1<b_1< a_2< \cdots <b_{k-1}< a_k <b_k.$$
We now wish to study the random vectors
\[\left(N_{g,n, \left[\frac{\sqrt{g}a_1}{n},\frac{\sqrt{g}b_1}{n}\right]}(X),\, \cdots,\, N_{g,n, \left[\frac{\sqrt{g}a_k}{n},\frac{\sqrt{g}b_k}{n}\right]}(X)\right).\]
For a tuple $\tuple{r_1}{r_k}\in\N$, define
\[A_{g,n} := \prod_{i=1}^{k}\left[\frac{\sqrt{g}a_i}{n},\frac{\sqrt{g}b_i}{n}\right]^{r_i}.\]

Now remark that because of the renormalisation of the lengths and that $\sqrt{g}/n\to 0$ as $g\to\infty$, each of the $\sqrt{g}b_i/n<2\mathrm{arcsinh}(1)$ for $g$ large enough. By the collar lemma, all closed geodesics shorter than $2\mathrm{arcsinh}(1)$ are disjoint and simple and so $Y_{g,n,r_1,\cdots,r_k}'=0$ for $g$ large enough. In particular, the proof of Theorem \ref{th: stat de poisson dans le cas n^2/h -> +infty} reduces to understanding only $\EE Y_{g,n,r_1,\cdots,r_k}$ for all $r_1,\ldots, r_k\in\N$.

The renormalisation of the lengths also requires a leading order expansion of the integrals
\[\int_{A_{g,n}} \prod_{i=1}^{K}\frac{4\sinh\left(\frac{x_i}{2}\right)^2}{x_i}\bigwedge_{i=1}^{K}\, \df x_i,\]
where $K=\sum_{i=1}^k r_i$ as usual.

\begin{lemma}
    \label{lemma: approximation integrale quand n^2/g -> infinity}
    For $r_1,\, \cdots, \, r_k\in\N$ and $K = \sum_{i=1}^k r_i$ we have
    \[\int_{A_{g,n}} \prod_{i=1}^{K}\frac{4\sinh\left(\frac{x_i}{2}\right)^2}{x_i}\bigwedge_{i=1}^{K}\, \df x_i =\left(\frac{g}{n^2}\right)^K\prod_{i=1}^{k}\left(\frac{b_i^2-a_i^2}{2}\right)^{r_i} + O\! \left(\left(\frac{g}{n^2}\right)^{K+1}\right).\]
\end{lemma}
\begin{proof}
    Note that,
    \[\int_{A_{g,n}} \prod_{i=1}^{k}\frac{4\sinh\left(\frac{x_i}{2}\right)^2}{x_i}\bigwedge_{i=1}^{k}\, \df x_i 
     = \prod_{i=1}^{k}\left(\int_{\frac{\sqrt{g}a_i}{n}}^{\frac{\sqrt{g}b_i}{n}}\frac{4\sinh\left(\frac{x}{2}\right)^2}{x}\, \df x\right)^{r_i}.\]
     By the substitution $u = (n/\sqrt{g})x$, we have
     \[\int_{\frac{\sqrt{g}a}{n}}^{\frac{\sqrt{g}b}{n}}\frac{4\sinh\left(\frac{x}{2}\right)^2}{x}\, \df x = \int_{a}^{b}\frac{4\sinh\left(\frac{\sqrt{g}u}{2n}\right)^2}{u}\, \df u.\]
     A first order Taylor approximation then furnishes the result
     \begin{align*}
        \int_{a}^{b}\frac{4\sinh\left(\frac{\sqrt{g}u}{2n}\right)^2}{u}\, \df u
        &=\int_{a}^{b}\frac{g}{n^2}u + O\!\left(\frac{g^2}{n^4}u^3\right)\, \df u
        =\frac{g}{n^2}\left(\frac{b^2-a^2}{2}\right) + O\! \left(\left(\frac{g}{n^2}\right)^{2}\right).
     \end{align*}
\end{proof}
We can now prove Theorem \ref{th: stat de poisson dans le cas n^2/h -> +infty}. The proof follows the same steps as the proof of Lemma \ref{lemma: n^2/g -> 0, multicurve} but this time using our new volume estimates. We will see that the main contribution is from the full pants multicurves.
\begin{proof}[Proof of Theorem \ref{th: stat de poisson dans le cas n^2/h -> +infty}]
    We will show 
    \[\EE\left(Y_{g,n,r_1\cdots, r_k}\right) = \prod_{i=1}^{k}\left(\frac{b_i^2 -a_i^2}{32\pi^2}\right)^{r_i} + O\!\left(\frac{1}{n}+ \frac{g}{n^2} + \frac{n}{g}\right).\]
    By Mirzakhani's integration formula Lemma \ref{lemma: Mirzakhani integration formula} we have
    \[\EE\left(Y_{g,n,r_1, \cdots, r_k}\right) = \sum_{q=1}^{K+1} \sum_{\substack{[\Gamma], \, \text{s.t.}\\q(\Gamma) = q}} 2^{-M(\Gamma)-K}\int_{A_{g,n}}\frac{V_{g,n}(\Gamma,x)}{V_{g,n}}x_1\cdots x_K \df x_1\wedge \cdots \wedge \df x_K.\]
    
    We first treat the case $q=1$. As in the proof of Lemma \ref{lemma: n^2/g -> 0, multicurve} we are dealing with the single mapping class group orbit of non-separating curves and so by Lemma \ref{Lemma 3.2 Mirzakhani2013} and the trivial $\sinh$ bound for bordered moduli space volumes we have
        $$\frac{V_{g,n}(\Gamma,x)}{V_{g,n}}=\frac{V_{g-K.n+2K}(x_1,x_1,\ldots,x_K,x_K)}{V_{g-K,n+2K}}\frac{V_{g-K,n+2K}}{V_{g,n}}\leq \prod_{i=1}^K\frac{4\sinh^2\left(\frac{x_i}{2}\right)}{x_i^2}.$$
    Thus by Lemma \ref{lemma: approximation integrale quand n^2/g -> infinity}, the $q=1$ term is bounded by 
    \begin{equation}
        \label{eq : cas n^2/g-> infinity q egale a 1}
        \sum_{\substack{[\Gamma], \, \text{s.t.}\\q(\Gamma) = 1}} 2^{-M(\Gamma)-K}\int_{A_{g,n}}\frac{V_{g,n}(\Gamma,x)}{V_{g,n}}x_1\cdots x_K \df x_1\wedge \cdots \wedge \df x_K\leq C\frac{g^{K}}{n^{2K}}.
    \end{equation} 
    For the terms $q=2,\ldots K$ we proceed as in the proof of Lemma \ref{lemma: n^2/g -> 0, multicurve} when dealing with the $q\neq 1$ terms, but this time we use the expansion for the $\sinh$ integral given by Lemma \ref{lemma: approximation integrale quand n^2/g -> infinity}. This gives 
    \begin{align}
    \label{eq : cas n^2/g-> infinity q entre 2 et K}
    &\sum_{q=2}^{K} \sum_{\substack{[\Gamma], \, \text{s.t.}\\q(\Gamma) = q}} 2^{-M(\Gamma)-K}\int_{A_{g,n}}\frac{V_{g,n}(\Gamma,x)}{V_{g,n}}x_1\cdots x_K \df x_1\wedge \cdots \wedge \df x_K\nonumber \\
    \leq& C\frac{g^K}{n^{2K}}\sum_{q=2}^{K}\sum_{\substack{a_1 + \cdots +a_q = 2K\\ a_i \geq 1}}
          \sum_{\substack{ I_1\sqcup\ldots\sqcup I_q=\{1,\ldots,n\}\\ g_1 + \cdots + g_q = g+q -K-1\\ 2g_i + n_i +a_i \geq \max(3, a_i)}}
        \frac{\prod_{i=1}^{q}V_{g_i,|I_i| + a_i}}{V_{g,n}}
          \prod_{i=1}^{K}\left(\frac{b_i^2 -a_i^2}{2}\right)^{r_i}\nonumber\\
          \leq& C\sum_{q=2}^{K}\sum_{\substack{a_1 + \cdots +a_q = 2K\\ a_i \geq 1}}\frac{g^{K-q+1}}{n^{2(k-q+1)}}\leq C\frac{g}{n^2},
    \end{align}
    where from the second to third lines we used Theorem \ref{th: somme d'un produit quelconque de V_gini}.
    
    We now deal with the term $q=K+1$ which contains the leading order contribution. Indeed, they will arise precisely from the contribution of full pants multicurves of length $K$. Let $\Gamma_{\mathrm{pants},K,K}$ denote the collection of mapping class group orbits of these multicurves. 
    
    Since each orbit in $\Gamma_{\mathrm{pants},K,K}$ provides the same contribution (since they all cut the same way, up to cusp labels), we can use Lemma \ref{lemma: nombre de pants multicurve} and Lemma \ref{lemma: Mirzakhani integration formula} to show that their total contribution is 
        $$\frac{n!}{2^K(n-2K)!}\int_{A_{g,n}}\frac{V_{g,n-K}\tuple{x_1}{x_K}}{V_{g,n}}\prod_{i=1}^{K}V_{0,3}(x_i)x_i\,\bigwedge_{i=1}^K \df x_i.$$
We know that $V_{0,3}(x_i)=1$. However, since $x_i\leq C\frac{\sqrt{g}}{n}$ we have
        $$1=\frac{2\sinh\left(\frac{x_i}{2}\right)}{x_i}+O\left(\frac{g}{n^2}\right).$$
We also have 
\begin{equation}
    \label{eq: approximation coeff binomiaux}
\frac{n!}{2^K(n-2K)!} = \frac{1}{2^K}n^{2K}\left(1+ O\!\left(\frac{1}{n}\right)\right),
\end{equation}
So, by using Lemma \ref{lemma: dl de Vgn(x)/Vgn} to expand the remaining volume, we see that this contribution is equal to 
    \begin{align*}
        \frac{n^{2K}}{2^K}\frac{V_{g,n-K}}{V_{g,n}}\int_{A_{g,n}}\prod_{i=1}^{k}\frac{4\sinh\left(\frac{x_i}{2}\right)^2}{x_i}\prod_{i=1}^{k}\, \df x_i + O\!\left(\frac{1}{n} +\frac{g}{n^2}+ \frac{n}{g}\right).
    \end{align*}
Next, by Theorem \ref{th: developpement premier ordre comparaison classique} we have
    $$\frac{V_{g,n-K}}{V_{g,n}}=\frac{1}{(8\pi^2g)^K}\left(1+O\left(\frac{n}{g}\right)\right).$$
Combining this with Lemma \ref{lemma: approximation integrale quand n^2/g -> infinity} provides
    \begin{align}
    \label{eq:leading-order n^2/g->infty}
    \frac{1}{(16\pi^2)^K}\prod_{i=1}^k\left(\frac{b_i^2-a_i^2}{2}\right)^{r_i}\left(1+O\left(\frac{1}{n}+\frac{n}{g}+\frac{g}{n^2}\right)\right),\end{align}
as desired.

It remains for us to show that the remaining terms for $q=K+1$ are in the error term.

To this end, to any multicurve mapping class group orbit $\Gamma$ with $q(\Gamma)=K+1$ and $\Gamma\notin\Gamma_{\mathrm{pants},K,K}$ we can associate to it the data $\{(g_i,a_i,I_i)\}_{i=1}^{K+1}$ where
    \begin{enumerate}
        \item $\sum_{i=1}^{K+1} a_i = 2K$, $a_i\geq 1$,
        \item $\sum_{i=1}^{K+1} g_i = g$, $g_i\geq 0$
        \item $I_1\sqcup\ldots\sqcup I_{K+1}=\{1,\ldots,n\}$,
        \item Some $m=0,\ldots,K$ for which $2g_i+|I_i|+a_i=3$ for $i=0,\ldots,m$ and $2g_i+|I_i|+a_i\geq \max\{4,a_i\}$ for $i=m+1,\ldots,K+1$,
        \item If $m=K$, then for all indices, we have $a_i\leq K-1$.
    \end{enumerate}
Indeed, this is just the separated surface data obtained after cutting $\Sigma_{g,n}$ by a multicurve in the orbit of $\Gamma$ similar to the proof of Lemma \ref{lemma: n^2/g -> 0, multicurve}. The first three conditions are obvious, the fourth condition tracks the number of components, $m$, that are pairs of pants. Note that $m\neq K+1$ because otherwise the surface would have Euler characteristic $-(K+1)$ which is never true when $g$ is sufficiently large. 

The fifth condition is forced due to the following. Obviously, every $a_i\leq K$ from the first condition. So, if there are $K$ pairs of pants among the $K+1$ components and there is some index satisfying $a_i=K$, then all $K$ of the pants components have $1$ boundary and $2$ cusps. But this means the multicurve is a full pants multicurve of length $K$ which we have excluded.

To any data we have (by identical reasoning as in the proof of Lemma \ref{lemma: n^2/g -> 0, multicurve}) at most a bounded (in terms of $K$) number of multicurve mapping class group orbits. 

By Mirzakhani's integration formula Lemma \ref{lemma: Mirzakhani integration formula}, the trivial $\sinh$ bound on bordered moduli space volumes and the expansion of the integral Lemma \ref{lemma: approximation integrale quand n^2/g -> infinity}, we see that the contribution of the remaining terms is bounded (up to some multiplicative constant depending only on $K$) by
\begin{align*}
    \frac{g^K}{n^{2K}}\sum_{m=0}^{K-1}\sum_{\substack{a_1 +\cdots+a_{K+1} = 2K\\ a_i \geq 1 \\ I_1\sqcup\ldots\sqcup I_{K+1}=\{1,\ldots,n\}\\ g_1 + \cdots + g_{K+1} = g\\ 2g_i +|I_i| +a_i = 3,\  i=0,\ldots,m\\
    2g_i+|I_i|+a_i\geq \max\{4,a_i\},\ i=m+1,\ldots,K+1}}\frac{\prod_{i=1}^{K+1}V_{g_i,|I_i| + a_i}}{V_{g,n}}
         \\
         +\frac{g^K}{n^{2K}}\sum_{\substack{a_1 +\cdots+a_{K+1} = 2K\\ 1\leq a_i\leq K-1 \\ I_1\sqcup\ldots\sqcup I_{K+1}=\{1,\ldots,n\}\\ g_1 + \cdots + g_{K+1} = g\\ 2g_i +|I_i| +a_i = 3,\  i=0,\ldots,K\\
    2g_{K+1}+|I_{K+1}|+a_{K+1}\geq \max\{4,a_i\}}}\frac{\prod_{i=1}^{K+1}V_{g_i,|I_i| + a_i}}{V_{g,n}}
\end{align*}
We now apply Theorem \ref{th: somme d'un produit quelconque de V_gini} to each of the inner summations (for fixed choices of the $a_i$, the number of choices of these $a_i$ is then bounded in terms of $K$ only) to obtain 
    \begin{align*}
        &\frac{g^K}{n^{2K}}\sum_{m=0}^{K-1}\sum_{\substack{a_1 +\cdots+a_{K+1} = 2K\\ a_i \geq 1 }}\sum_{i=1}^{K+1}\frac{1+n^{\sum_{\substack{j=1\\j\neq i}}^m\max\{a_j,3\}-a_j+\sum_{\substack{j=m+1\\j\neq i}}^{K+1}\max\{a_j,4\}-a_j}}{g^{-a_i+\sum_{\substack{j=1\\j\neq i}}^m\max\{a_j,3\}-a_j+\sum_{\substack{j=m+1\\j\neq i}}^{K+1}\max\{a_j,4\}-a_j}}\\
        \\
         &\hspace{1cm}+\frac{g^K}{n^{2K}}\sum_{\substack{a_1 +\cdots+a_{K+1} = 2K\\ 1\leq a_i\leq K-1 }}\sum_{i=1}^{K+1}\frac{1+n^{\sum_{\substack{j=1\\j\neq i}}^K\max\{a_j,3\}-a_j+\max\{a_{K+1},4\}-a_{K+1}}}{g^{-a_i+\sum_{\substack{j=1\\j\neq i}}^K\max\{a_j,3\}-a_j+\max\{a_{K+1},4\}-a_{K+1}}}\\
         &=\frac{g^{3K}}{n^{4K}}\sum_{m=0}^{K-1}\sum_{\substack{a_1 +\cdots+a_{K+1} = 2K\\ a_i \geq 1 }}\sum_{i=1}^{K+1}\frac{1+n^{a_i+\sum_{\substack{j=1\\j\neq i}}^m\max\{a_j,3\}+\sum_{\substack{j=m+1\\j\neq i}}^{K+1}\max\{a_j,4\}}}{g^{\sum_{\substack{j=1\\j\neq i}}^m\max\{a_j,3\}+\sum_{\substack{j=m+1\\j\neq i}}^{K+1}\max\{a_j,4\}}}\\
        \\
         &\hspace{1cm}+\frac{g^{3K}}{n^{4K}}\sum_{\substack{a_1 +\cdots+a_{K+1} = 2K\\ 1\leq a_i\leq K-1 }}\sum_{i=1}^{K+1}\frac{1+n^{a_i+\sum_{\substack{j=1\\j\neq i}}^K\max\{a_j,3\}+\max\{a_{K+1},4\}}}{g^{\sum_{\substack{j=1\\j\neq i}}^K\max\{a_j,3\}+\max\{a_{K+1},4\}}},
    \end{align*}
where the equality uses the condition $\sum_{i=1}^Ka_i=2K$. Now notice that the summands are maximised by taking $a_i$ maximal and the other $a_j$ minimal. In the first term this corresponds to taking $a_i=K$ and the other $a_j=1$ and assuming that the $i$th term is always in the second summation in the exponents of $n$ and $g$; in the second term it corresponds to $a_i=K-1$ and the other $a_j=1$. Up to multiplicative constant depending only on $K$ we result in 
     \begin{align}
     \label{eq: q=K+1-error}
         &\frac{g^{3K}}{n^{3K}}\sum_{m=0}^{K-1}\sum_{\substack{a_1 +\cdots+a_{K+1} = 2K\\ a_i \geq 1 }}\frac{1+n^{4K-m}}{g^{4K-m}}+\frac{g^{3K}}{n^{3K+1}}\sum_{\substack{a_1 +\cdots+a_{K+1} = 2K\\ 1\leq a_i\leq K-1 }}\frac{1+n^{3K}}{g^{3K}}\leq C\left(\frac{n}{g}+\frac{1}{n}\right),
    \end{align}
which indeed lies in the error term. The result now follows from combining \eqref{eq : cas n^2/g-> infinity q egale a 1}, \eqref{eq : cas n^2/g-> infinity q entre 2 et K}, \eqref{eq:leading-order n^2/g->infty} and \eqref{eq: q=K+1-error}.
\end{proof}

\subsection{Proof of Theorem \ref{thm:stat de poisson dans le cas n^2/g -> alpha}}

In this section we will prove Theorem \ref{thm:stat de poisson dans le cas n^2/g -> alpha} which works under the assumption that
    $$\lim_{g\to\infty}\frac{n^2}{g}\to\alpha\in (0,\infty).$$
Let $k,r_1,\ldots,r_k\in\N$ and
    $$0\leq a_1<b_1< \cdots < a_k < b_k < 2\mathrm{arcsinh}(1).$$
Set $K = \sum_{i=1}^k r_i$ as usual and
\[\alpha_i = \frac{\alpha}{4\pi^2}\left(\cosh\left(\frac{b_i}{2}\right)-\cosh\left(\frac{a_i}{2}\right)\right) \qquad \text{and} \qquad \lambda_i = \int_{a_i}^{b_i} \frac{2\sinh^2\left(\frac{x}{2}\right)}{x} \, \df x.\]
Note that this time we will only be considering closed geodesics with lengths less than $2\mathrm{arcsinh}(1)$ to ensure that they are simple and disjoint. This is because, in this regime of $n$ and $g$, non-simple curves will also contribute to the leading order. Indeed, this can be seen in the proof of Lemma \ref{lemma: self-insection negligeable} where the error term is of order $n^2/g$ and it is not hard to see that there are terms of this order (for example those coming from non-simple closed curves that fill a pair of pants).

For these shorter closed geodesics, we will still recover Poisson statistics in terms of the parameters $\alpha_i$ and $\beta_i$.

\begin{proof}[Proof of Theorem \ref{thm:stat de poisson dans le cas n^2/g -> alpha}]

As remarked, since we consider closed geodesics of length at most $2 \mathrm{arcsinh} (1)$, we have \[(N_{g,n,[a_1,b_1]})_{r_1}\cdots(N_{g,n,[a_k,b_k]})_{r_k} = Y_{g,n,r_1, \ldots, r_k}.\] 
Using Mirzakhani's integration formula Lemma \ref{lemma: Mirzakhani integration formula} we deduce that 
\begin{align*}
    \EE\left(Y_{g,n,r_1, \cdots, r_k}\right) = \sum_{q=1}^{K+1} \sum_{\substack{[\Gamma], \, \text{s.t.}\\q(\Gamma) = q}} 2^{-M(\Gamma)-K}\int_{\prod_{i=1}^{k}[a_i,b_i]^{r_i}}\frac{V_{g,n}(\Gamma,x)}{V_{g,n}}x_1\cdots x_K \df x_1\wedge \cdots \wedge \df x_K.
\end{align*}
As usual, we consider first the case when $q=1$. This follows identically to the proof of Lemma \ref{lemma: n^2/g -> 0, multicurve} because the error term is of order $n/g$ and we have no scaling in the lengths. So with applications of Lemma \ref{lemma: dl de Vgn(x)/Vgn} and Theorem \ref{th: developpement premier ordre comparaison classique} we once again obtain 
\begin{align}
\label{eq:n^2=g q=1}
&2^{-M(\Gamma_0)-K}\int_{\prod_{i=1}^{k}[a_i,b_i]^{r_i}}\frac{V_{g,n}(\Gamma_0,x)}{V_{g,n}}\left(\prod_{i=1}^{K}x_i\right) \bigwedge_{i=1}^{K} \df x_i\nonumber \\
 =& \frac{1}{2^K}\frac{V_{g-K,n+2K}}{V_{g,n}}\int_{\prod_{i=1}^{k}[a_i,b_i]^{r_i}}\frac{V_{g-K,n+2K}(x_1,x_1, \cdots, x_K,x_K)}{V_{g-K,n+2K}}\left(\prod_{i=1}^{K}x_i\right) \bigwedge_{i=1}^{K} \df x_i\nonumber\\
 =&\prod_{i=1}^{k}\lambda_i^{r_i}\left(1 + O\! \left(\frac{n}{g}\right)\right).
\end{align}

The leading order contribution now arises in all of the $q=2,\ldots,K+1$ terms from the mapping class group orbits of $(q-1)$-pants multicurves of length $K$. Recall from Lemma \ref{lemma: nombre d'orbite de pants} that once we fix the positions of the $q-1$ separating curves in such a multicurve orbit, there are precisely 
    $$\frac{n!}{(n-(2(q-1))!}=n^{2(q-1)}\left(1+O\left(\frac{1}{n}\right)\right),$$
mapping class group orbits of these multicurves. Moreover, after cutting the pants off, the remaining component has genus $g+q-K-1$, $2K-(q-1)$ boundaries and $n-2(q-1)$ cusps.

So, now we apply Mirzakhani's integration formula Lemma \ref{lemma: Mirzakhani integration formula} to this class of multicurve orbits to obtain a contribution of
    \begin{align*}
        \sum_{q=2}^{K+1}\sum_{1\leq i_1<i_2<\ldots<i_{q-1}\leq K}\frac{n^{2(q-1)}}{2^K}\int_A\prod_{i=1}^Kx_i\frac{V_{g+q-K-1,n-3(q-1)+2K}(\mathbf{x}_{i_1,\ldots,i_{q-1}},\tilde{\mathbf{x}}_{i_1,\ldots,i_{q-1}},\tilde{\mathbf{x}}_{i_1,\ldots,i_{q-1}})}{V_{g,n}}\bigwedge_{i=1}^K\df x_i,
    \end{align*}
where the inner summation ranges over the choices of positions for the separating components of the multicurve and we use the notation
    $$\mathbf{x}_{i_1,\ldots,i_{q-1}}=(x_{i_1},\ldots,x_{i_{q-1}})\ \text{and}\ \tilde{\mathbf{x}}_{i_1,\ldots,i_{q-1}}=(x_j)_{j\neq i_\ell,\ \forall \ell=1,\ldots,q_1}.$$
Note that we have also used the fact that $V_{0,3}(x,y,z)=1$. This differs from the proof of Theorem \ref{th: stat de poisson dans le cas n^2/h -> +infty} where we are able to approximate these volumes by the usual $2\sinh(x/2)/x$ up to small errors because the $x$ were sufficiently small to be in the error term.

By Lemma \ref{lemma: dl de Vgn(x)/Vgn} and Theorem \ref{th: developpement premier ordre comparaison classique}, we have 
    \begin{align*}
        &\frac{V_{g+q-K-1,n-3(q-1)+2K}(\mathbf{x}_{i_1,\ldots,i_{q-1}},\tilde{\mathbf{x}}_{i_1,\ldots,i_{q-1}},\tilde{\mathbf{x}}_{i_1,\ldots,i_{q-1}})}{V_{g,n}}\\
        &=\frac{V_{g+q-K-1,n-3(q-1)+2K}}{V_{g,n}}\prod_{j=i_1,\ldots,i_{q-1}}\frac{2\sinh\left(\frac{x_j}{2}\right)}{x_j}\prod_{j\neq i_1,\ldots,i_{q-1}}\frac{4\sinh^2\left(\frac{x_j}{2}\right)}{x_j^2}\left(1+O\left(\frac{n}{g}\right)\right)\\
        &=\frac{1}{(8\pi^2g)^{q-1}}\prod_{j=i_1,\ldots,i_{q-1}}\frac{2\sinh\left(\frac{x_j}{2}\right)}{x_j}\prod_{j\neq i_1,\ldots,i_{q-1}}\frac{4\sinh^2\left(\frac{x_j}{2}\right)}{x_j^2}\left(1+O\left(\frac{n}{g}\right)\right).
    \end{align*}
In particular, the contribution becomes
    \begin{align}
    \label{eq:n^2=g leading q}
        &\sum_{q=2}^{K+1}\sum_{1\leq i_1<\ldots<i_{q-1}\leq K}\left(\prod_{j=1}^{q-1}\frac{n^2}{8\pi^2g}\int_{[a_{i_j},b_{i_j}]}\sinh\left(\frac{x}{2}\right)\df x\right)\left(\prod_{\substack{i=1 \\ i\neq i_1,\ldots,i_{q-1}}}^K\int_{[a_i,b_i]}\frac{2\sinh^2\left(\frac{x}{2}\right)}{x}\df x\right)\left(1+O\left(\frac{n}{g}\right)\right)\nonumber\\
        &=\sum_{q=1}^{K}\sum_{1\leq i_1<\ldots<i_{q}\leq K}\left(\prod_{j=1}^{q}\frac{n^2}{4\pi^2g}\left(\cosh\left(\frac{b_{i_j}}{2}\right)-\cosh\left(\frac{a_{i_j}}{2}\right)\right)\right)\left(\prod_{\substack{i=1 \\ i\neq i_1,\ldots,i_{q}}}^K\lambda_i\right)\left(1+O\left(\frac{n}{g}\right)\right).
    \end{align}
Notice that \eqref{eq:n^2=g q=1} corresponds precisely to $q=0$ in the summation in \eqref{eq:n^2=g leading q}. So combining these terms, we recognise an expansion of the form
    $$\prod_{i=1}^K(A_i+B_i)=\sum_{q=0}^K\sum_{1\leq i_1<\ldots<i_q\leq K}\prod_{j=1}^q A_{i_j}\prod_{\substack{j=1 \\ j\neq i_1,\ldots,i_q}}^K B_j.$$
Letting
    $$A_i(n(g),g) = \frac{n^2}{4\pi^2g}\left(\cosh\left(\frac{b_i}{2}\right)-\cosh\left(\frac{a_i}{2}\right)\right),$$
we obtain that the combined contributions are
    \begin{align}
        \prod_{j=1}^K\left(A_j(n(g),g\right)+\lambda_j)\left(1+O\left(\frac{n}{g}\right)\right)=\prod_{j=1}^k\left(A_j(n(g),g\right)+\lambda_j)^{r_j}\left(1+O\left(\frac{n}{g}\right)\right).
    \end{align}

The remaining multicurve orbits that we must consider are those that are not $(q-1)$-pants multicurves of length $K$. However, we can use an identical argument as in the proof of Theorem \ref{th: stat de poisson dans le cas n^2/h -> +infty} onwards from equation \eqref{eq:leading-order n^2/g->infty} (after multiplying through by $n^{2K}/g^K$ since we no longer rescale the geodesic lengths) to show that they have combined contribution that is 
    $$O\left(\frac{n}{g}+\frac{1}{n}\right),$$
since none of the error terms are worse than this and this remains unchanged after the multiplication by $n^{2K}/g^K$ and using that $n^2\asymp g$.

Thus, since $A_i(n(g),g)\to\alpha_i$ as $g\to\infty$, we obtain that
\begin{equation*}
\EE\left((N_{g,n,[a_1,b_1]})_{r_1}\cdots(N_{g,n,[a_k,b_k]})_{r_k}\right) \to  \prod_{i=1}^{k}(\alpha_i+\lambda_i)^{r_i},
\end{equation*}
as $g\to\infty$ and the result follows from Theorem \ref{th : Method of factorial moments}.
\end{proof}
Although we don't demonstrate Poisson statistics here when the lengths of the geodesics are longer than $2\mathrm{arcsinh}(1)$ we immediately obtain in any case that they have finite factorial moments of all orders.

\begin{corollary}
\label{cor: finite-moments}
    Let $k\in\N$ and $0\leq a_1<b_1< \ldots < a_k < b_k$. Then, for any $r_1, \ldots, r_k\in\N$, 
    \[\lim_{g\to \infty} \EE\left((N_{g,n,[a_1,b_1]})_{r_1} \cdots (N_{g,n,[a_k,b_k]})_{r_k}\right) < \infty.\]
\end{corollary}
\begin{proof}
    Because $\EE\left((N_{g,n,[a_1,b_1]})_{r_1} \cdots (N_{g,n,[a_k,b_k]})_{r_k}\right) = Y_{g,n, r_1, \cdots, r_k} + Y_{g,n,r_1, \cdots, r_k}^{\prime}$, the proofs of Theorem \ref{thm:stat de poisson dans le cas n^2/g -> alpha} and Lemma \ref{lemma: self-insection negligeable} give the result.
\end{proof}

\section{Systoles and small eigenvalues}
In this section, we prove Proposition \ref{prop: systole2} and Theorem \ref{thm: eigenavlues}. The proof of Corollary \ref{cor: systole} is identical to \cite[Proposition 5.1]{MirzakhaniPetri2021}. For both of these results, we are working in the regime
    \[\lim_{g \to \infty}\frac{n(g)}{g} = 0 \qquad \text{and} \qquad \lim_{g \to \infty}n(g) = \infty.\]
For $a<b$, we introduce the random variable
    $$N_{g,n,2,[a,b]}:\mathcal{M}_{g,n}\to\Z_{\geq 0},$$
that counts the number of primitive, unoriented and simple closed geodesics with length in $[a,b]$ that separate a pair of pants with two cusps from the surface. In other words, it counts the number of full pants multicurves of length $1$ with lengths in $[a,b]$.

\begin{proof}[Proof of Proposition \ref{prop: systole2}]
The proof follows the lines of \cite[Proposition 1.4]{HideThomas2024} or \cite[Appendix A]{MirzakhaniPetri2021}. Since the systole is always a simple closed geodesic, for any $x>0$ we have
    $$\PP\left(\mathrm{systole}(X)>\frac{t\sqrt{g}}{n}\right)\leq\PP\left(N_{g,n,2,\left[0,\frac{t\sqrt{g}}{n}\right]}(X)=0\right).$$
Since the sequence $x_g$ that we consider is $o\left(\frac{n}{\sqrt{g}}\right)$, we are only considering curves that are length at most $2\mathrm{arcsinh}(1)$ which are disjoint by the collar lemma. For any $k\geq 1$, we can thus compute the factorial moments identically to the contribution of full pants multicurves as in the proof of Theorem \ref{th: stat de poisson dans le cas n^2/h -> +infty} to be
    \begin{align*}
        \EE\left(\left(N_{g,n,2,\left[0,\frac{t\sqrt{g}}{n}\right]}(X))\right)_k\right) &=\frac{1}{2^k}\frac{1}{V_{g,n}}{n\choose \underbrace{2,\ldots, 2}_{\text{$k$ times}}}\int_{\left[0,\frac{t\sqrt{g}}{n}\right]^k} V_{g,n-k}(x_1,\ldots,x_k)\prod_{i=1}^k x_i\bigwedge_{i=1}^k\df x_i\\
        &=\frac{n!}{2^k(n-2k)!}\frac{V_{g,n-k}}{V_{g,n}}\int_{\left[0,\frac{t\sqrt{g}}{n}\right]^k} \prod_{i=1}^k \sinh\left(\frac{x_i}{2}\right)\left(1+O\left(\frac{t^2}{n}\right)\right)\bigwedge_{i=1}^k\df x_i,        
    \end{align*}
where on the second line we have used Lemma \ref{lemma: dl de Vgn(x)/Vgn}. By Theorem \ref{th: developpement premier ordre comparaison classique}, 
    $$\frac{V_{g,n-k}}{V_{g,n}}= \left(\frac{1}{8\pi^2g}\right)^k\left(1+O\left(\frac{n}{g}\right)\right),$$
and     
    $$\int_{0}^{\frac{t\sqrt{g}}{n}}\sinh\left(\frac{x}{2}\right)\df x= 2\left(\cosh\left(\frac{t\sqrt{g}}{2n}\right)-1\right)=\frac{t^2g}{2n^2}\left(1+O\left(\frac{t^2g}{n^2}\right)\right).$$
So put together,
    $$\EE\left(\left(N_{g,n,2,\left[0,\frac{t\sqrt{g}}{n}\right]}(X))\right)_k\right) = \frac{t^{2k}}{(32\pi^2)^k}\left(1+O\left(t^2\left(\frac{1}{n}+\frac{g}{n^2}\right)+\frac{n}{g}\right)\right).$$
Let $\lambda_g=\frac{x_g^2}{32\pi^2}$. By Markov's inequality, we have
    \begin{align*}
        \PP&\left(N_{g,n,2,\left[0,\frac{x_g\sqrt{g}}{n}\right]}(X))=0\right)\\
        &\leq \PP\left(\left|N_{g,n,2,\left[0,\frac{x_g\sqrt{g}}{n}\right]}(X))-\lambda_g\right|\geq \lambda_g\right)\\
        &\leq \frac{\EE\left(N_{g,n,2,\left[0,\frac{x_g\sqrt{g}}{n}\right]}(X))-\lambda_g\right)^2}{\lambda_g^2}\\
        &=\frac{1}{\lambda_g^2}\left(\EE\left(\left(N_{g,n,2,\left[0,\frac{x_g\sqrt{g}}{n}\right]}(X))\right)_2\right)  +(1-2\lambda_g)\EE\left(N_{g,n,2,\left[0,\frac{x_g\sqrt{g}}{n}\right]}(X)\right)+\lambda_g^2\right)\\
        &=\frac{32\pi^2}{x_g^2}+O\left(x_g^2\left(\frac{1}{n}+\frac{g}{n^2}\right)+\frac{n}{g}\right),
    \end{align*}
as desired.
\end{proof}

We now turn to the proof of Theorem \ref{thm: eigenavlues} which follows the proof of \cite[Theorem 1.1]{Hi.Th2025} with the appropriate modifications. We recall the following result on showing the existence of small eigenvalues.

\begin{lemma}[{\cite[Lemma 5.3]{HideThomas2024}}]
\label{lem: curves-to-evalues}
Let $\varepsilon<\frac{1}{4}$ and assume that $X\in\mathcal{M}_{g,n}$ has $k+1$ simple closed geodesics of length at most $\varepsilon/6$, each separating off a pair of pants with two cusps. Then $X$ has least $k$ eigenvalues and the $k$th largest eigenvalue satisfies $\lambda_k\leq \varepsilon$.
\end{lemma}

\begin{proof}[Proof of Theorem \ref{thm: eigenavlues}]
    We will show that given $0<L<2\mathrm{arcsinh}(1)$, there exists $C(L)>0$ such that
        $$\PP\left(\frac{gN_{g,n,2,\left[0,L\right]}(X))}{n^2}<C(L)\right)\to 0,$$
    as $g\to\infty$. Under this assumption, if we let $L=\varepsilon/6$, then 
        $$N_{g,n,2,[0,\frac{\varepsilon}{6}]}(X)>\frac{n^2}{g}C(\tfrac{\varepsilon}{6})$$
    with probability tending to $1$ as $g\to\infty$ and so by Lemma \ref{lem: curves-to-evalues}, we obtain
        $$\lambda_{\frac{n^2}{g}C(\tfrac{\varepsilon}{6})}<\varepsilon,$$
    with probability tending to $1$ as $g\to\infty$.
    
    To prove the claim, we can appeal to the same computation as in the proof of Proposition \ref{prop: systole2}. We have
        $$\EE\left(\frac{gN_{g,n,2,[0,L]}}{n^2}\right)=\frac{1}{8\pi^2}\left(\cosh\left(\frac{L}{2}\right)-1\right)\left(1+O\left(\frac{n}{g}+\frac{1}{n}\right)\right).$$
    Similarly,
        \begin{align*}
    \mathrm{Var}\left(\frac{gN_{g,n,2,[0,L]}}{n^2}\right)&= \frac{g^2}{n^4}\EE\left(\left(N_{g,n,2,[0,L]}(X)\right)_2\right)+\frac{g^2}{n^4}\EE(N_{g,n,[0,L]}(X))-\frac{g^2}{n^4}\EE\left(N_{g,n,[0,L]}(X)\right)^2\\
    &=O\left(\frac{n}{g}+\frac{1}{n}\right)
        \end{align*}
    Now, take $C'(L)<\frac{1}{8\pi^2}\left(\cosh\left(\frac{L}{2}\right)-1\right)$. Then
    Chebyshev's inequality gives
        $$\PP\left(\left|\frac{gN_{g,n,2,[0,L]}(X)}{n^2}-\EE\left(\frac{gN_{g,n,2,[0,L]}(X)}{n^2}\right)\right|\geq C'(L)\right)= O\left(\frac{n}{g}+\frac{1}{n}\right).$$
    So, with probability at least $1-O\left(\frac{n}{g}+\frac{1}{n}\right)$ as $g\to\infty$, we have
        $$\frac{gN_{g,n,2,[0,L]}(X)}{n^2}\geq \EE\left(\frac{gN_{g,n,2,[0,L]}(X)}{n^2}\right)-C'(L)\geq C(L),$$
    for some $0<C(L)<\frac{1}{8\pi^2}\left(\cosh\left(\frac{L}{2}\right)-1\right)-C'(L)$.
\end{proof}

\appendix
\section{Proof of Lemma \ref{lemma 5.1 MZ}}
\label{Appendix: proof lemma 5.1}
For Mirzakhani and Zograf, the proof of Lemma \ref{lemma 5.1 MZ} essentially follows from an estimate as in Theorem 3.3(ii) in their article \cite{MirzakhaniZograf2024}. In this proof however, they use as input an estimate on sums of products of volumes akin to our Lemma \ref{lemme initialisation de rec produit des volumes} that is valid for fixed $n$ (see \cite[Lemma 3.3]{Mirzakhani2013}) and so does not include explicit dependency on $n$ in the errors. Inputting the bound from Lemma \ref{lemme initialisation de rec produit des volumes} instead remedies this and requires no further modification to the schematic of their proof. We include the details for completeness but stress that the ideas are all contained in the paper of Mirzakhani and Zograf \cite{MirzakhaniZograf2024}.

The starting point is an easy estimate on differences of intersection numbers (actually only the case where $n_0=1$ is needed).
\begin{lemma}
    \label{lemma: Th 3.2 (ii) MZ}
    Let $n_0 \in \N$ such that $n_0 \leq n$, let $\tuple{d_1}{d_{n_0}} \in \N^{n_0}$ and assume that
\[\lim_{g\to \infty} \frac{n}{g} = 0.\] 
Then, there exists a universal constant $C>0$ for which 
\begin{align*}
\left[ \tau_{d_1} \cdots \tau_{d_{n_0}} \tau_0^{\,n-n_0} \right]_{g,n}
- \left[ \tau_{d_1 +1} \cdots \tau_{d_{n_0}} \tau_0^{\,n-n_0} \right]_{g,n}
&\leq  C\left(\left(\sum_{j =1}^{n_0} d_j + nd_1\right)V_{g,n-1} +d_1V_{g-1, n+1}\right).
\end{align*}
\end{lemma}
\begin{proof}
    Recall, as in the proof of Lemma \ref{lemma expansion 1 er ordre dif intersection number}, that we have the decomposition
   \begin{align*}
   \left[\tau_{d_1}\cdots\tau_{d_{n_0}}\tau_{0}^{n-n_0}\right]_{g,n} - \left[\tau_{d_1+1}\cdots\tau_{d_{n_0}}\tau_{0}^{n-n_0}\right]_{g,n}
   =& \tilde{A}_{d,g,n} + \tilde{B}_{d, g,n} + \tilde{C}_{d,g,n},
   \end{align*}
   where $\tilde{A}_{d,g,n},\, \tilde{B}_{d, g,n}\ \text{and} \ \tilde{C}_{d,g,n}$ are written explicitly in the proof therein.
   
   The terms $\tilde{A}_{d,g,n}$ and $\tilde{B}_{d, g,n}$ are bounded in exactly the same way as in the proof of Theorem 3.3 of \cite{MirzakhaniZograf2024}. First we trivially bound the intersection numbers by their corresponding volumes and then bound the sum of the $a_l$ differences using Lemma \ref{lemme:  aproxiamtion des coefficient a_i} to obtain 
   \[\tilde{A}_{d,g,n} \leq 8 \left(2\sum_{j =2}^{n_0} d_j + n-1\right)V_{g, n-1}\]
   and 
   \[\tilde{B}_{d, g,n} \leq 16d_1V_{g-1, n+1}.\]
   We proceed similarly for $\tilde{C}_{d,g,n}$ by bounding the intersection numbers by volumes but then apply Lemma \ref{lemme initialisation de rec produit des volumes} to obtain 
   \[\dsum_{\substack{I\sqcup J = \ens{2}{n}\\ g_1 +g_2 = g}} V_{g_1, \, \left|I\right| +1}V_{g_2, \, \left|J\right| +1}\leq C\frac{n^2}{g}.\]
   So,
   \[\tilde{C}_{d,g,n} \leq C\sum_{i=0}^{\infty}(a_{i+1}-a_i)(i+d_1+1)\frac{n^2}{g} V_{g,n-1}\leq C'nd_1V_{g,n-1}.\]
   And because $n/g$ goes to $0$ as $g$ goes to infinity and thanks to \ref{lemme: aproxiamtion des coefficient a_i},
   \[\tilde{C}_{d,g,n} \leq  4(1+2d_1)nV_{g,n-1}.\]
   Adding the three terms and we get the result.
\end{proof}
We now prove Lemma \ref{lemma 5.1 MZ}.
\begin{proof}
    We prove (i) in a more general case as promised in Remark \ref{remark: on lemma 5.1 MZ}. Let $n_0, d_1, \, \cdots, \, d_{n_0} \in \N$. Then, as in the proof of Theorem \ref{th: developpement premier ordre comparaison classique} we have
    
    \begin{align*}
1 - \frac{\left[\tau_{d_1}\cdots\tau_{d_{n_0}}\tau_{0}^{n-n_0}\right]_{g,n}}{V_{g,n}}    =& \underbrace{\sum_{k=0}^{d_1-1}\frac{\left[\tau_k\tau_0^{n-1}\right]_{g,n}-\left[\tau_{k+1}\tau_0^{n-1}\right]_{g,n}}{V_{g,n}}}_{(1)}\\
      &+\underbrace{\sum_{i=1}^{n_0-1}\sum_{k=0}^{d_{i+1}-1}\frac{\left[\tau_{k}\tau_{d_1}\cdots\tau_{d_i}\tau_{0}^{n-(i+1)}\right]_{g,n}-\left[\tau_{k+1}\tau_{d_1}\cdots\tau_{d_i}\tau_{0}^{n-(i+1) }\right]_{g,n}}{V_{g,n}}}_{(2)}.
\end{align*}
    So, thanks to Lemma \ref{lemma: Th 3.2 (ii) MZ}, we bound (1) by  
    $$\frac{C}{V_{g,n}}\left(nd_1^2V_{g,n-1}+d_1^2V_{g-1,n+1}\right)\leq C|d|^2\frac{nV_{g,n-1}+V_{g-1,n+1}}{V_{g,n}},$$
    and (2) by
    \begin{align*}
        \frac{C}{V_{g,n}}\sum_{i=1}^{n_0-1}\sum_{k=0}^{d_{i+1}-1}\left(|d|nV_{g,n-1}+d_{i+1}V_{g-1,n+1}\right)&\leq \frac{C}{V_{g,n}}|d|\left(nV_{g,n-1}+V_{g-1,n+1}\right)\sum_{i=1}^{n_0-1} d_{i+1}\\
        &=C|d|^2\left(\frac{nV_{g,n-1}+V_{g-1,n+1}}{V_{g,n}}\right).
    \end{align*}
    Put together and then using Lemma \ref{Lemma 3.2 Mirzakhani2013} then gives the desired upper bound.
    
    The lower bound is very coarse, using the bound
        $$\left[\tau_{d_1}\cdots\tau_{d_{n_0}}\tau_0^{n-n_0}\right]_{g,n}\leq \left[\tau_1\tau_0^{n-1}\right]_{g,n},$$ 
        along with the recursion Lemma \ref{lemma: recursion_rels}\eqref{eq: rec III} and Lemma \ref{Lemma 3.2 Mirzakhani2013} to obtain 
    \begin{align*}
      \left(1 - \frac{\left[\tau_{d_1}\cdots\tau_{d_q}\tau_{0}^{n-q}\right]_{g,n}}{V_{g,n}}\right)
      &\geq \left(1 - \frac{\left[\tau_{1}\tau_{0}^{n-1}\right]_{g,n}}{V_{g,n}}\right)\\
      &\geq \frac{\tilde{A}_{0,g,n}}{V_{g,n}}\\
      &\geq 8(n-1)(a_1 -a_0)\frac{V_{g,n-1}}{V_{g,n}}\\
      &\geq 8(a_1-a_0)b_0\frac{n-1}{2g+n-3}\\
      &\geq C\frac{n-1}{2g+n-2}. 
   \end{align*}
   So for $n \geq 2$ we get the result.  For (ii), recall that
   \[\sum_{l=1}^{+\infty}\frac{l\pi^{2l-2}}{(2l+1)!} = \frac{\pi \cosh(\pi) - \sinh(\pi)}{2\pi^3}.\]
   So by Lemma \ref{lemma: recursion_rels}\eqref{eq: rec II} and the upper bound proven above, we have 
   \[\left|\frac{(2g+n-2)V_{g,n}}{V_{g,n+1}}- \frac{1}{4\pi^2}\right| \leq \frac{1}{2}\sum_{l=1}^{+\infty}\frac{l\pi^{2l-2}}{(2l+1)!}\left|\frac{\left[\tau_{l}\tau_{0}^{n}\right]_{g,n+1}}{V_{g,n+1}}-1\right| \leq c_2\frac{\max(1,n)}{2g+n-2}.\]
   Finally, for (iii) we use the recursion Lemma \ref{lemma: recursion_rels}\eqref{eq: rec I} along with the lower bound obtained in the proof of (i) to obtain
      \begin{align*}
      \frac{V_{g-1,n+4}}{V_{g,n+2}}\leq \frac{\left[\tau_{1}\tau_{0}^{n+1}\right]_{g,n+2}}{V_{g,n+2}}\leq 1- c_0\frac{n+1}{2g+n} \leq 1- c_3\frac{n}{2g+n-2}.
   \end{align*}
This concludes the proof.
\end{proof}

\bibliographystyle{abbrv}
\bibliography{biblio.bib }

\begin{tabular}{@{}l@{}}
Malo Gatignol \\
CentraleSupélec,
Université Paris-Saclay,
Gif-sur-Yvette, France\\
\texttt{malo.gatignol@student-cs.fr}
\end{tabular}

\vspace{1em}

\begin{tabular}{@{}l@{}}
Joe Thomas \\
Department of Mathematical Sciences,
Durham University, 
Durham, UK \\
\texttt{joe.thomas@durham.ac.uk}
\end{tabular}
\end{document}